\documentclass[12pt,journal,onecolumn]{IEEEtran} 

\title{\bf
Convergence of Nonlinear Observers on $\RR ^n$
with a Riemannian Metric (Part II)
}
\author{Ricardo G. Sanfelice\thanks{R. G. Sanfelice is with
the 
Department of Computer and Engineering, University of California,
1156 High Street, Santa Cruz, CA 95064.
       Email: {\tt\small ricardo@ucsc.edu}. Research partially supported by NSF Grant no. ECS-1150306 and by AFOSR Grant no. FA9550-12-1-0366.
}
and
Laurent Praly\thanks{L. Praly is with CAS, ParisTech, Ecole des Mines, 35 rue Saint Honor\'{e}, 77305, 
        Fontainebleau, France
       Email: {\tt\small Laurent.Praly@ensmp.fr}
}
}

\usepackage{amsmath}
\usepackage{latexsym}
\usepackage{amsfonts}
\usepackage[dvipsnames]{xcolor}
\usepackage{psfrag}
\usepackage{graphics}
\usepackage{graphicx}
\usepackage{setspace}
\usepackage{hyperref}
\usepackage{ifthen}
\usepackage{subfigure}
\usepackage{hyperref}
\usepackage{url}
\usepackage{balance}

\DeclareMathAlphabet\EuScript{U}{eus}{m}{n}
\SetMathAlphabet\EuScript{bold}{U}{eus}{b}{n}

\def\euA{{\EuScript{A}}}
\def\euB{{\EuScript{B}}}
\def\euC{{\EuScript{C}}}
\def\euH{{\EuScript{H}}}
\def\euKL{\EuScript{K}_{\nu}}

\def\euPL{\EuScript{P}\!_{\nu}}
\def\euP{\EuScript{P}}

\def\gr{\textrm{grad}}
\def\Hess{\textrm{Hess}}

\def\dummyone{\mathfrak{b}}
\def\dummytwo{\mathfrak{a}}
\def\dummythree{\mathfrak{c}}
\def\dummyfour{\mathfrak{d}}

\def\addots{\mathinner{\mskip1mu\raise 1pt\vbox{\kern 7pt\hbox{.}}\mskip2mu
    \raise 4pt\hbox{.}\mskip2mu\raise 7pt\hbox{.}\mskip1mu}}

\newcommand{\ox}{\xi}

\def\XX{X}
\newcommand{\reals}{\mathbb{R}}

\def\RR{\mathbb{R}}
\def\SS{\mathbb{S}}

\def\Ouv{\Omega }
\newcommand{\A}{{\cal{A}}}
\newcommand{\bigbrace}[1]{\left\{#1\right\}}
\newcommand{\defset}[2]{\bigbrace{#1\ : \ #2 }}
\newcommand{\x}{x}
\newcommand{\xdot}{\dot{\x}}
\newcommand{\f}{f}
\newcommand{\h}{h}
\newcommand{\y}{y}
\newcommand{\F}{F}

\newcommand{\xhat}{\hat{x}}
\newcommand{\xhatdot}{\dot{\xhat}}

\catcode`\@=11
\def\downparenfill{$\m@th\braceld\leaders\vrule\hfill\bracerd$}
\def\overparen#1{\mathop{\vbox{\ialign{##\crcr\crcr \noalign{\kern0.4ex}
\downparenfill\crcr\noalign{\kern0.4ex\nointerlineskip}
$\hfil\displaystyle{#1}\hfil$\crcr}}}\limits}
\catcode`\@=12

\newcommand{\RRgeq}{[0,+\infty)}

\newcommand{\plower}{\underline{p}}

\newcommand{\pupper}{\overline{p}}
\newcommand{\qlower}{q}

\def\Did{{\mathfrak{D}^+}}

\def\xrond{{\mathchoice%
{{\mbox{$\scriptstyle \mathcal{X}$}}}%
{{\mbox{$\scriptstyle \mathcal{X}$}}}%
{{\mbox{$\scriptscriptstyle \mathcal{X}$}}}%
{{\mbox{\tiny$\scriptscriptstyle \mathcal{X}$}}}%
}}

\newtheorem{theorem}{Theorem}[section] 
\newtheorem{lemma}[theorem]{Lemma}   
\newtheorem{example}[theorem]{Example}    
\newtheorem{proposition}[theorem]{Proposition}    
\newtheorem{remark}[theorem]{Remark}   
\newtheorem{definition}[theorem]{Definition}

\begin{document}

\maketitle

\begin{abstract}
In \cite{Sanfelice.Praly.11},
it is established that a convergent observer with an infinite
gain margin can be designed for a given nonlinear system
when a Riemannian metric showing that the system is differentially detectable 
(i.e.,~the Lie derivative of the Riemannian
metric along the system vector
field is negative in the space tangent to
the output function level sets) and 
the level sets of the output function are geodesically convex is available.
In this paper, we propose techniques for designing a Riemannian metric 
satisfying the first property in the case where the system is strongly infinitesimally observable
(i.e., each time-varying linear system
resulting from the linearization along a
solution to the system satisfies a uniform observability property) 
or where it is strongly differentially observable (i.e. the mapping state to output 
derivatives is an injective immersion) or where it is Lagrangian.
Also, we give results 
that are complementary to those in \cite{Sanfelice.Praly.11}. In particular,
we provide a locally convergent observer and  make a link to the existence of
a reduced order observer. Examples illustrating the results are presented.
\end{abstract}

\section{Introduction}

We consider a nonlinear system of the form
\footnote{%
\label{footnote1}
If the system is time varying (perhaps due to known exogenous inputs), i.e.,
$\dot x = f(x,t)$, $y = h(x,t)$
most of the results of \cite{Sanfelice.Praly.11} as well as those here can be 
extended readily by simply replacing $x$ by 
$x_e = [x^\top\ t]^\top$, leading to the time-invariant system with dynamics
$
\dot x_e = [f(x,t)^\top\ 1]^\top =: f_e(x_e)$, $y_e = [h(x,t)^\top\ t]^\top =: h_e(x_e).$
The drawback of this simplifying viewpoint is that, when time 
dependence is induced by exogenous inputs, for each input we 
obtain a different time-varying system.
And, maybe even more handicapping, we need to know the time-variations 
for the design.
}
\begin{equation}
\label{eqn:Plant1}
\xdot  \;=\;  \f(\x)\  ,\quad 
\y  \;=\;  \h(\x),
\end{equation}
with $\x $ in $\RR ^n$ being the system's state and $\y$ in $\RR ^p$ the 
measured system's output. We are interested in the design of a 
function $\F$ such that
the set
\begin{equation}\label{eqn:A}
\A := \defset{(\x, \xhat)\in\RR ^n\times\RR ^{n}}{ x = \xhat}
\end{equation}
is asymptotically stable
for
the system
\begin{equation}
\label{LP56}
\dot x \;=\;  f(x)\quad ,\qquad 
\xhatdot \;=\;   \F(\xhat,h(x))
\  .
\end{equation}

A solution to this problem that was proposed in \cite{Sanfelice.Praly.11}
is re-stated in Theorem \ref{thm1part1},
which is in Section~\ref{sec:LieCondition}.
It relies on the formalism of Riemannian geometry 
and gives conditions under which a constructive procedure exists for 
getting an appropriate function $\F$.
This solution requires the satisfaction of mainly two conditions. 
The first condition is about the geodesic convexity 
 of the level sets of the output function
 (see point \ref{LP73} in Appendix~\ref{sec:Glossary}). This condition
is not addressed here.
Instead,
 we focus our attention on the second condition,
which is a 
differential detectability property\footnote{This expression was suggested
to us by Vincent Andrieu.}, made precise in Definition \ref{def:DifferentialDetectability}
below. With the terminology used in the study of contracting 
flows in Riemannian spaces, this property means that
$f$ is strictly geodesically monotonic tangentially to the output function 
level sets. 
Forthcoming examples related to the so-called harmonic oscillator
with unknown frequency will illustrate these notions and provide
metrics certifying both weak and strong differential detectability.

In Section \ref{sec:LieCondition}, we establish results complementing 
those in \cite{Sanfelice.Praly.11}.
In Section~\ref{sec:LocalConvergenceCase}, we establish that
the differential detectability property only is already sufficient
to obtain a locally convergent observer.
In Section \ref{sec1}, we show that this property implies also the existence of a locally  convergent
reduced order observer, in this way, extending
the result established in  \cite[Corollary 3.1]{Besancon.00.SCL} for the 
particular case where the metric is Euclidean.
The conclusion we draw from Section \ref{sec:LieCondition} is that the design of a locally 
convergent observer can be reduced  to the design a 
metric exhibiting the differential stability property.
Sections \ref{sec:Reconstructibility}, \ref{sec:High gain}, and 
\ref{ex:Lagrangian} are dedicated to such designs in three different
contexts.

In Section~\ref{sec:Reconstructibility},
under a uniform observability property of
the family of time-varying linear systems resulting from the linearization along solutions to the system, 
a symmetric covariant $2$-tensor giving the strong 
differential detectability property is shown to exist as a solution 
to a Riccati equation
which, for linear systems, would be an algebraic Riccati equation.
Proposition~\ref{propo:NecessitySatisfactionReconstructibilityDRE} establishes this fact.
The resulting metric leads to an observer that
resembles the Extended Kalman Filter; 
see, e.g., \cite{jazwinski1970stochastic}.
In Section~\ref{sec:Reconstructibility},
Proposition~\ref{prop:PidotLambda} shows that the metric can instead be taken in
the form of an exponentially weighted observability Grammian, 
leading to an observer design method that is in the spirit of the one proposed in 
\cite{Kleinman.70}. 

In Section~\ref{sec:High gain},
for systems that are strongly differentially observable
\cite[Chapter 2.4]{Gauthier.Kupka.01.Book.Observers}, 
we propose an expression for the tensor
that is based on the fact that, after writing the system dynamics in an observer form, a high 
gain observer can be used.
This result leads to an observer which has some similarity with the 
observer for linear systems obtained using  Ackerman's formula.

Finally, in Section~\ref{ex:Lagrangian},
we show how a Riemannian metric can be constructed for Euler-Lagrange 
systems whose Lagrangian is quadratic in the generalized velocities.
This result extends the result in \cite{Aghannan.Rouchon.03}.

The design methods proposed in Section~\ref{sec:Reconstructibility} do not 
necessarily lead to explicit expressions for the metric. Instead, they give 
numerical procedures to compute it, only involving the solution of
ordinary differential equations over a grid of initial conditions. 
On the other hand,
the designs in Sections~\ref{sec:High gain} and \ref{ex:Lagrangian} involve computations that can be done symbolically. 
All of these various designs are
coordinate independent and do not require to have the
system written in some specific form.

To ease the reading, we give a glossary in Appendix~\ref{sec:Glossary}
definitions of the main objects we employ from differential geometry.

\section{Full and Reduced Observers under Strong Differential Detectability}
\label{sec:LieCondition}
In this section, we study what can be obtained when the system 
satisfies the differential detectability  property defined as follows
(see items \ref{LP65}, \ref{LP73}, and \ref{LP63} in Appendix~\ref{sec:Glossary}).
\begin{definition}
\label{def:DifferentialDetectability}
The nonlinear system \eqref{eqn:Plant1} is {\em strongly differentially detectable} (respectively, {\em weakly differentially detectable})
on a closed,
weakly geodesically convex set $\mathcal{C} \subset \RR^n$
with nonempty interior
if there exists
a symmetric covariant $2$-tensor $P$ on $\RR^n$ satisfying 
\begin{equation}
\label{eqn:LfPtangentially}\renewcommand{\arraystretch}{1.4}
\begin{array}{ll}\displaystyle
v^\top \mathcal{L}_fP(x) v
<0 \qquad (\mbox{respectively } \leq 0) & \\ 
& \hspace{-1.8in} \displaystyle \qquad \forall (x,v) \in \mathcal{C}\times\SS^{n-1}\,  :\: dh(x)v=0\ .
\end{array}
\end{equation}
\end{definition}

We illustrate this property with an example

\begin{example}
\label{ex0}
Consider a harmonic oscillator with unknown frequency. Its 
dynamics are
\begin{equation}
\label{eqn:HarmonicOscillator}
\dot x\;=\;  f(x)\; :=\; 
\left(\renewcommand{\arraystretch}{0.8}\begin{array}{c}
x_2 \\ -x_3\,  x_1 \\ 0
\end{array}\right)
\quad ,\qquad 
 y\;=\; h(x) := x_1
\end{equation}
with $(x_1,x_2,x_3) \in \reals\times\reals\times \RR_{>0}$.
As a candidate to check the differential detectability we pick, in 
the above coordinates,
\begin{equation}
\label{LP1}
P(x) = \left(
\begin{array}{ccc}
1 + 2 \ell k^2 + 4 \ell^2 x_1^2 & 
- 2 \ell k & 2 \ell x_1  \\ 
- 2 \ell k & 2 \ell & 0 \\
2 \ell x_1  & 0 & 1
\end{array}\right)
\  .
\end{equation}
where $k$ and $\ell$ are strictly positive real numbers.
The expression of its Lie derivative $\mathcal{L}_f P$ in these 
coordinates is
\\[1em]\vbox{\noindent%
\null \hfill $\displaystyle 
\left(\begin{array}{ccc}
4\ell k x_3 + 8 \ell^2 x_1 x_2
&
\star
&
\star
\\
1 + 2 \ell k^2 + 4 \ell^2 x_1^2-2\ell x_3
&
-4\ell k
&
\star
\\
2\ell k x_1 + 2 \ell x_2 
&
0
&
0
\end{array}\right)
$\hfill\null}\\[1em]
where
the various $\star$ should be replaced by their symmetric 
values.
Then, since we have
$
\frac{\partial h}{\partial x}(x)v = v_1
$, where $v = (v_1,v_2,v_3)$,
the evaluation of the Lie derivative of $P$  for a vector $v$ in the kernel 
of $dh$ gives
\begin{equation}
\label{LP82}
\left(\begin{array}{cc}
v_2 & v_3
\end{array}\right)
\left(\begin{array}{ccc}
-4\ell k
&
0
\\
0
&
0
\end{array}\right)
\left(\begin{array}{@{}c@{}}
v_2 \\ v_3
\end{array}\right)
\;=\; -4\ell k v_2^2
\  .
\end{equation}
This allows us to conclude that the harmonic oscillator with unknown 
frequency is weakly differentially detectable. Actually,
as we shall see later when we use a different metric, 
it is strongly differentially detectable.
\hfill $\triangle$
\end{example}

With this property of differential detectability at hand, we study in 
the next two subsections what it implies in terms of existence of 
converging full and then reduced order observers.

\subsection{Local Asymptotic Stabilization of the set $\A$}
\label{sec:LocalConvergenceCase}
In \cite[Theorem 3.3 and Lemma 3.6]{Sanfelice.Praly.11} we have established the following result
(see also \cite{Sanfelice.Praly.13.ArXiv}).

\begin{theorem}
\label{thm1part1}
Assume there exist a Riemannian metric $P$ and
a closed subset $\mathcal{C} $ of $\RR ^n$, with 
nonempty interior,
such that
\begin{enumerate}
\item[A1~:]
$\mathcal{C}$ 
is weakly geodesically convex;
\item[A2~:]
There exist 
a continuous function $\rho :\RR ^n\to \RRgeq$
and a strictly positive real number $q$ such that
\begin{equation}
\label{3}
\mathcal{L}_f P(\x)
\; \leq  \; 
\rho (\x)\,  
d\h (x)\otimes d\h(\x) 
\;-\; \qlower\,  P(\x)
\qquad \forall \x\in\mathcal{C}
\  ;
\end{equation}
\item[A3~:]
There exists
a $C^2$ function 
$\RR^p\times\RR^p \ni  (y_a,y_b) \mapsto\delta(y_a,y_b) \in \RRgeq$
satisfying
$$
\delta (h(x),h(x))\;=\; 0 , \ \
\left.
\frac{\partial ^2\delta }{\partial y_a^2}(y_a,y_b)
\right|_{y_a=y_b=h(x)}\; >\; 0
$$
for all $x\in \mathcal{C}$,
and, such that, for any pair $(x_a,x_b)$ in $\mathcal{C}\times 
\mathcal{C}$ satisfying
$
h(x_a) \neq h(x_b)
$
and, for any minimizing geodesic $\gamma ^*$ between
$x_a=\gamma ^*(s_a)$ and $x_b=\gamma ^*(s_b)$
satisfying
$\gamma ^*(s) \in \mathcal{C}$ for all $s $ in $ [s_a, s_b]$,
$s_a\leq s_b$,
we have
$$
\frac{d}{ds}\delta (h(\gamma ^*(s)),h(\gamma ^*(s_a)))\; >\; 0
\qquad \forall s\in (s_a,s_b]
\  .
$$
\end{enumerate}
Then,
 for any positive real number $E$ there exists 
a continuous function $k_E:\RR ^n\to \RR $ such that,
with the observer given by (see item 
\ref{LP66} in Appendix \ref{sec:Glossary})
\begin{equation}\label{eqn:GeodesicObserverVectorField}
F(\xhat,y)\;=\; f(\xhat) \;-\; k_E(\xhat)\,  \gr_P h(\xhat)
\frac{\partial \delta }{\partial y_a}(\h(\xhat),y)^\top\  ,
\end{equation}
the following holds\footnote{%
$\Did d(\xhat,\x)$ is the upper right Dini derivative along the 
solution, i.e., with $(\hat X((\hat x,\x),t),X(x,t))$ denoting a 
solution of (\ref{LP56}),
$$
\Did d(\xhat,\x)\;=\; 
\limsup_{t\searrow 0} \frac{d(\hat X((\hat x,\x),t),X(x,t))-d(\hat x,x)}{t}
$$
}:
$$
\Did d(\xhat,\x)
\; \leq \;  \displaystyle -
\frac{\qlower}{4} 
\,  d(\xhat,\x)$$ 
for all $(x,\xhat) \in
\left\{(x,\xhat )\,  :\:
d(\xhat,x) <  E
\right\}
\; \bigcap\; 
\left(\mbox{\rm int}(\mathcal{C}) \times \mbox{\rm int}(\mathcal{C})\right).
$
\end{theorem} 

Theorem~\ref{thm1part1} establishes that, when assumptions A1-A3 hold,
for every given positive number $E$,
an observer with vector field as in \eqref{eqn:GeodesicObserverVectorField}
renders the set $\A$ in \eqref{eqn:A} asymptotically stable with
a domain of attraction containing the set
$$\left\{(x,\xhat )\,  :\: d(\xhat,x) <  E \right\} \; \bigcap\; 
\left(\mbox{\rm int}(\mathcal{C}) \times \mbox{\rm int}(\mathcal{C})\right)$$

Condition A2 is a stronger version of what we have called 
differential detectability in the introduction. We come back to it 
extensively below.

Condition A3 is a restrictive way of saying that the output level sets 
are geodesically convex.
Fortunately, even without assumption A3, 
inspired by \cite[Theorem 1]{Aghannan.Rouchon.03},
we can design 
an observer making the set \eqref{eqn:A} asymptotically stable.
As opposed to Theorem~\ref{thm1part1},
its domain of attraction cannot be made arbitrarily large.

\begin{proposition}
\label{prop1}
Assume there exist a Riemannian metric $P$ and
a closed subset $\mathcal{C} $ of $\RR ^n$, with 
nonempty interior,
such that
\begin{enumerate}
\item[A1'~:]
$\mathcal{C}$ 
is weakly geodesically convex and there exist coordinates denoted 
$\x$ and positive numbers
 $\plower$ and $\bar h_1$
such that, for each $\x  $ in $ \mathcal{C}$, we have
\begin{equation}\label{eqn:BoundsForLocalResult}
\begin{array}{c}
\displaystyle 
\plower\; \leq \; |P(\x )|
\quad ,
\qquad  
\left|\Hess _P h(x)\right|\; \leq \; \bar h_1\
\end{array}
\end{equation}
where $\Hess _P h$ is the $p$-uplet of the Hessian
of the components $h_i$ of $h$;
see item~\ref{item:CovariantDerivativeOfh} in Appendix~\ref{sec:Glossary}.
\item[A2'~:]
There exist a positive real number $\bar \rho $
and a strictly positive real number $\underline{q}$ such that
\begin{equation}
\label{3withConstantRho}
\mathcal{L}_f P(\x )
\; \leq  \; 
\bar \rho \,    
dh(x)\otimes dh(x)
\;-\; \underline{q} \,  P(\x )
\qquad \forall \x \in\mathcal{C}.
\end{equation}
\item[A3'~:]
There exists
a $C^2$ function
$\RR^p\times\RR^p \ni (y_a,y_b) \mapsto\delta(y_a,y_b) \in \RRgeq$
and positive real numbers $\bar \delta _1$ and 
$\underline{\delta }_2$
satisfying
\begin{eqnarray}
\label{LP94}
\displaystyle \delta (h(x),h(x))\;=\; 0, \ \ 
\left.
\frac{\partial ^2\delta }{\partial y_a^2}(y_a,y_b)
\right|_{y_a=y_b=h(x)}\;\hspace{-0.3in} >\; \underline{\delta }_2 \,  I \hspace{-0.2in}
\end{eqnarray}
for all $x\in\mathcal{C}$,
\begin{eqnarray}
\displaystyle\label{LP95}
\left|\frac{\partial \delta }{\partial y_a}(h(x_a),h(x_b))\right|
\; \leq \; \bar \delta _1\,  d(x_a,x_b)
\end{eqnarray}
for all $(x_a,x_b)\in \mathcal{C}\times\mathcal{C}$.
\end{enumerate}
Then, 
with the observer given by
\begin{equation}
\label{4}
F({\hat{x}},y)\;=\; f({\hat{x}}) \;-\; k\,  \gr_Ph(\xhat)
\frac{\partial \delta }{\partial y_a}(h(\hat x),y)^\top
\  ,
\end{equation}
the following holds:
\begin{equation}
\label{5}
\Did d({\hat{x}},x )
\; \leq \; 
-\underline{r}
\,  d({\hat{x}},x ) 
\end{equation}
for all $(x,{\hat{x}}) \in
\left\{(x,{\hat{x}} )\,  :\:
d({\hat{x}},x) \leq   \frac{\varepsilon }{k}
\right\}
\; \bigcap\; 
(\mathcal{C} \times \mathcal{C})$
when we have
\begin{equation}
\label{LP67}
k\geq \frac{\bar \rho }{2 \underline{\delta }_2}
\quad ,\qquad 
\underline{q} > \underline{r} 
\quad ,\qquad 
\varepsilon \;:=\; 
\frac{(\underline{q}-\underline{r})\plower }{ 2\bar h_1 \bar\delta_1}
\  .
\end{equation}
\end{proposition}

\begin{remark}
We make the following observations:
\begin{enumerate}
\item
A key difference with respect to the result in Theorem~\ref{thm1part1} is that,
in the latter, the domain of attraction gets larger with the increase of the observer gain,
while the domain of attraction guaranteed by the result in Proposition~\ref{prop1}
decreases when $k$ increases.
\item
When there exists a positive real number $\bar h_2$ satisfying
$$
\left|\frac{\partial h}{\partial \x }(\x )\right|\; \leq \; \bar h_2 
\qquad \forall x \in \mathcal{C}
\  ,
$$
a function $\delta $ satisfying A3' is simply
$$
\delta (y_a,y_b)\;=\; |y_a-y_b|^2
$$
Indeed, let $\gamma ^*:[s_a,s_b] \to \reals^n$ be a minimizing 
geodesic between $x_a$ and $x_b$ that stays in $\mathcal{C}$. 
We have
\\[0.7em]
\null\vbox{\hsize=\linewidth\noindent
$
\begin{array}{@{}cl@{}}\displaystyle
\left|\frac{\partial \delta }{\partial y_a}(h(x_a),h(x_b))\right| =& 
\displaystyle
2 \left|h(x_a)-h(x_b)\right|
\  ,
\\
\hspace{1.2in}=&
\displaystyle 
2\left|\int_{s_a}^{s_b}
\frac{\partial h}{\partial x}(\gamma ^*(r))\frac{d\gamma ^*}{ds}(r) dr\right|
\  ,
\\[1.5em]
\hspace{1.2in}=&\displaystyle 
2\int_{s_a}^{s_b}
\sqrt{\frac{\partial h}{\partial x}(\gamma ^*(r))
P(\gamma ^*(r))^{-1}
\frac{\partial h}{\partial x}(\gamma ^*(r))^\top
}\\[1.5em]
&\displaystyle 
\hspace{0.8in} \times
\sqrt{\frac{d\gamma ^*}{ds}(r) ^\top P(\gamma ^*(r))\frac{d\gamma 
^*}{ds}(r) }\: dr
\  ,
\\[0.7em]
\hspace{1.2in}\leq &\displaystyle 
\frac{2\bar h_2}{\sqrt{\plower }}\,
d(x_a,x_b)
\  .
\end{array}
$}\\[0.5em]
\end{enumerate}
\hfill $\Box$
\end{remark}

\begin{proof}
It is sufficient to show that the vector field $\hat x\mapsto F(\hat x,y)$ is 
geodesically strictly monotonic with respect to $P$ (uniformly in $y$), at least when $\hat x$ and $x$ are sufficiently close.  See \cite[Lemma 2.2]{Sanfelice.Praly.11} and the discussion before it.
With the coordinates given by assumption A1',
and item~\ref{item:CovariantDerivativeOfh} in 
Appendix~\ref{sec:Glossary},
we have
\begin{eqnarray*}
\mathcal{L}_FP(\hat x,y)
&\hskip -0.3em=&\hskip -0.3em
\mathcal{L}_f P(\hat x)
\:-\: k \mathcal{L}_{\gr _P h} P(\hat x,y) \otimes
\frac{\partial \delta }{\partial y_a}(h(\hat x),y)^\top 
 \\ 
 &\hskip -0.3em &\hskip -0.3em
 -\: 2k\,  \frac{\partial h}{\partial x}({\hat{x}})^\top 
\frac{\partial ^2\delta }{\partial y_a^2}(h(\hat x),y)\frac{\partial h}{\partial x}({\hat{x}}) 
\:  ,
\\
&\hskip -0.3em = &\hskip -0.3em
\mathcal{L}_f P(\hat x)
\:-\: 2k\,  \Hess _P h (\hat x) 
\otimes
\frac{\partial \delta }{\partial y_a}(h(\hat x),y)^\top
 \\ 
 &\hskip -0.3em &\hskip -0.3em
-\: 2k\,  \frac{\partial h}{\partial x}({\hat{x}})^\top 
\frac{\partial ^2\delta }{\partial y_a^2}(h(\hat x),y)\frac{\partial h}{\partial x}({\hat{x}}) 
\:  .
\end{eqnarray*}
Here, the notation $\Hess _P h \otimes v$, with $v$ a vector in $\RR^p$ stands for
$
\sum_{i=1}^p \Hess _P h_i\: v_i
\,  ,
$
where each $\Hess _P h_i\,   v_i$ is a covariant $2$-tensor.
So, with  \eqref{eqn:BoundsForLocalResult}, (\ref{3withConstantRho}), 
(\ref{LP94}), (\ref{LP95}) and (\ref{LP67}), we obtain successively
\begin{eqnarray*}
\mathcal{L}_FP(\hat x,y)& \leq & 
\mathcal{L}_f P(\hat x)
\;+\; 
 2k\,  \bar h_1 \bar \delta _1 d(\hat x,x) 
\;-\;  2k\,  \underline{\delta }_2\, 
\frac{\partial h}{\partial x}({\hat{x}})^\top \frac{\partial h}{\partial x}({\hat{x}}) 
\  ,
\\
&\leq &
-\underline{q} P(\hat x)
\;+\; k\,  \frac{2\bar h_1\bar\delta_1}{\plower }
d(\hat x,x)\, P(\hat x)
\;-\; (2k\underline{\delta }_2-\bar \rho )\, 
\frac{\partial h}{\partial x}({\hat{x}})^\top \frac{\partial h}{\partial x}({\hat{x}}) 
\  ,
\\
&\leq &
 -\underline{r} P(\hat x)
\end{eqnarray*}
for all $(x,{\hat{x}}) \in
\left\{(x,{\hat{x}} )\,  :\:
d({\hat{x}},x) \leq   \frac{\varepsilon }{k}
\right\}
\; \bigcap\; 
(\mathcal{C} \times \mathcal{C})$.
Since $\mathcal{C}$ is weakly geodesically convex, (\ref{5}) follows by 
integration along a minimizing geodesic.
\end{proof}
The proofs of Theorem \ref{thm1part1} and Proposition \ref{prop1} 
differ mainly on the way the term $\Hess _P h (\hat x) 
\otimes
\frac{\partial \delta }{\partial y_a}(h(\hat x),y)^\top$ is handled. 
With Assumption A3, related to the geodesic convexity 
of the output level sets, it can be shown to be harmless because of 
its sign. Instead, with Assumption A3' only, we go with upper bounds 
and show it is harmless at least when $\hat x$ and $x$ are 
sufficiently close. Hence, a local convergence result in the latter case and a 
regional  one in the former are obtained.

\subsection{A Link between the Existence of $P$ and a Reduced Order Observer}
\label{sec1}

In \cite[Corollary 3.1]{Besancon.00.SCL} it is established that, if, 
in some coordinates,
the expression of the metric $P$ is constant and that of $h$ is linear, 
then there exists a reduced order observer.  In this section, we 
establish a similar result without imposing the metric to be Euclidean.
The interest of a reduced order observer is that there is no 
correction term to design. This task is replaced by that of finding 
appropriate coordinates. In our context, the existence of such 
coordinates is guaranteed by
the following result from \cite{Eisenhart.25}.

\begin{theorem}[{\cite[p. 57 \S 19]{Eisenhart.25}}]
\label{thm:Eisenhart2}
Let $P$ be a complete Riemannian metric on $\RR^n$.
Assume $p=1$ and  $h$ has rank $1$ at $x_0 $ in $\RR^n$. Then, there exists a neighborhood 
$\mathcal{N}_{x_0}$ of 
$x_0$ on which there exists coordinates
$$
x\;=\; (y,\xrond )
$$
such that,
for each $x $ in $ \mathcal{N}_{x_0}$, the expression of $h$ and $P$ in 
these coordinates can be decomposed as
\begin{equation}
\label{EQU9}
y\;=\; h((y,\xrond ))
\end{equation}
and
\begin{equation}
\label{EQU8}
P((y,\xrond ))\;=\; \left(\begin{array}{cc}
P_{yy}(y,\xrond )& 0
\\
0
&
P_{\xrond\xrond}(y,\xrond )
\end{array}\right)
\  ,
\end{equation}
with
$P_{yy}(y,\xrond) $ in $\reals^{p \times p}$ and
$P_{\xrond\xrond}(y,\xrond)$ in $\reals^{(n-p) \times (n-p)}$.
\end{theorem}

\begin{proof}
See \cite[p. 57 \S 19]{Eisenhart.25}. A sketch
of another proof is as follows. Note first that,
the Constant Rank Theorem  implies the existence of
a neighborhood of $x_0$ on which coordinates $(y ,\bar \xrond )$
are defined and satisfy
$
h(x) \;=\; h((y ,\bar \xrond ))\;=\; y
\  .
$
Let the expression of the metric in the $(y,\bar \xrond)$-coordinates 
be
$$
\overline{P}((y,\bar \xrond  ))
=
\left(
\begin{array}{@{}cc@{}}
\overline{P}_{yy}(y,\bar \xrond  ) & \overline{P}_{y\bar \xrond  }(y,\bar \xrond  )
\\[0.5em]
\overline{P}_{\bar \xrond  y}(y,\bar \xrond  ) & \overline{P}_{\bar \xrond  \bar \xrond  }(y,\bar \xrond  )
\end{array}\right)
$$
and let $\varphi(y,\bar \xrond)$ denote the solution, evaluated at 
time $h(x_0)$, of the 
time-varying system
$
\frac{d\mathfrak{x}}{d \mathfrak{y}}\;=\;
-\overline{P}_{\bar \xrond\bar \xrond}(\mathfrak{y},\mathfrak{x})^{-1}
\overline{P}_{\bar \xrond y}(\mathfrak{y},\mathfrak{x})
$
issued from $\mathfrak{x}=\bar \xrond$ at time $\mathfrak{y}=y$.
The proof can be completed by showing that the function $\varphi$ 
defined this way on a neighborhood of 
$x_0$ satisfies all the required properties for
$
(y,\xrond)\;=\; (y\,  ,\,  \varphi(y ,\bar \xrond))
$
to be the appropriate coordinates
in the neighborhood of $x_0$.
\end{proof}
\begin{example}
\label{ex2}
Consider the matrix $P$ in (\ref{LP1}) with
$
y=x_1$, $\bar \xrond=(x_2,x_3)$.
We have
$$
\overline{P}_{\bar \xrond  y}(y,\bar \xrond  )=
\left(\begin{array}{@{}c@{}}
-2\ell k\\2\ell y
\end{array}\right)
\  ,\quad 
\overline{P}_{\bar \xrond  \bar \xrond  }(y,\bar \xrond  )
=\left(\begin{array}{cc}
2 \ell & 0 \\ 0 & 1
\end{array}\right)
$$
This leads to the system
$$
\frac{d\mathfrak{x}}{d \mathfrak{y}}\;=\;\mathfrak{f}\left(\mathfrak{y},\mathfrak{x}\right)\;=\; 
-\overline{P}_{\bar \xrond\bar \xrond}(\mathfrak{y},\mathfrak{x})^{-1}
\overline{P}_{\bar \xrond y}(\mathfrak{y},\mathfrak{x})
\;=\; \left(\begin{array}{@{}c@{}}
k \\ -2\ell \mathfrak{y}  
\end{array}\right)
$$
the solutions of which, at time $\mathfrak{y}$, going through $\mathfrak{x}_0$ at time 
$\mathfrak{y}_0$, are
$$
\mathfrak{X}(\mathfrak{x}_0,\mathfrak{y}  _0;\mathfrak{y} )
\;=\; \mathfrak{x}_0
\;+\; \left(\begin{array}{@{}c@{}}
 k [\mathfrak{y} -\mathfrak{y} _0]
\\
-\ell [\mathfrak{y} ^2-\mathfrak{y} _0^2]
\end{array}\right)
$$
So in particular, we get
$$
\varphi ((y ,\bar \xrond))
=
\mathfrak{X}((x_2,x_3),y;0 )
=
\left(\begin{array}{@{}c@{}}
x_2 - k y
\\
x_3 + \ell y^2
\end{array}\right)
\  .
$$
From the proof above, it follows that the coordinates 
$(y,\xrond)$ satisfying (\ref{EQU8}) in Theorem \ref{thm:Eisenhart2}
are
\begin{equation}
\label{LP2}
\left(
y, \xrond _1, \xrond _2
\right)
=
\varphi (x)
=
 \varphi ((y ,\bar \xrond))
= \left(
x_1, x_2 - k x_1, x_3 + \ell x_1^2
\right)
\  .
\end{equation}
They are defined on the open set
\begin{equation}
\label{LP7}
\Ouv\;=\; \mathcal{N}_{x_0}=\varphi (\RR^2\times\RR_{>0})
\end{equation}
and they give
\\[0.7em]\null \hfill $\displaystyle 
P_{yy}((y,\xrond))\;=\; 1
\  ,\quad 
P_{\xrond\xrond}((y,\xrond))\left(\begin{array}{cc}
2\ell & 0 \\ 0 & 1
\end{array}\right)\  .
$
\hfill $\triangle$
\end{example}

\par\vspace{1em}
Let us express the differential detectability and 
the observer (\ref{eqn:GeodesicObserverVectorField}) in the special 
coordinates given by Theorem \ref{thm:Eisenhart2}. The dynamics of \eqref{eqn:Plant1} in 
the coordinates $(y,\xrond )$ are
$$
\dot y = f_y(y,\xrond )
\quad ,\qquad 
\dot {\xrond } = f_\xrond (y,\xrond )
$$
We notice that, by decomposing a 
tangent vector as
$
v = \left(\begin{array}{@{}c@{}}
v_y \\ v_\xrond
\end{array}\right),
$
and since $\frac{\partial h}{\partial y}(x_0) \not = 0$,
we find that \eqref{EQU9} gives,
for every $x =(y,\xrond)$ in $\mathcal{N}_{x_0}$,
$$
\frac{\partial h}{\partial x}(x ) v = 0 
\ \ \
\Longleftrightarrow 
\ \ \ 
\frac{\partial h}{\partial y}(y,\xrond ) v_y = 0 
\ \ \
\Longleftrightarrow 
\ \ \
v_y = 0
\  .$$
It follows that, with expression (\ref{EQU8}) 
and in $(y,\xrond )$ coordinates,
condition A2 in (\ref{3}) 
is as follows:
\\[1em]$\displaystyle 
2\,  v_\xrond ^\top
P_{\xrond \xrond} (y,\xrond )
\frac{\partial f_\xrond}{\partial \xrond}(\xrond )v_\xrond 
+ 
\frac{\partial }{\partial y}
\left(
\vrule height 0.6em depth 0.6em width 0pt
v_\xrond ^\top P_{\xrond\xrond}(y,\xrond  )v_\xrond \right)
f_y(y,\xrond )
+
\frac{\partial }{\partial \xrond}
\left(
\vrule height 0.6em depth 0.6em width 0pt
v_\xrond ^\top P_{\xrond\xrond}(y,\xrond  )v_\xrond \right)
f_\xrond(y,\xrond)
\leq  \: - \qlower\,  v_\xrond ^\top P_{\xrond\xrond}(y,\xrond  )v_\xrond 
$\  \null\refstepcounter{equation}\label{LP57}$(\theequation)$\\[1em]
for all $(y,\xrond ,v_\xrond)$ such that $(y,\xrond)\in \mathcal{N}_{x_0}, \; v_\xrond \in  \SS^{n-2}$.
Also our observer  
(\ref{eqn:GeodesicObserverVectorField}) takes the form
\begin{eqnarray*}
\dot {\hat y }&=& f_{y}(\hat y,\hat {\xrond })
\;-\;
k_E((\hat y,\hat {\xrond }))
\frac{1}{P_{yy}((\hat y,\hat {\xrond }))}
\frac{\partial \delta }{\partial y_a}(\hat y ,y)
\  ,
\\[-0.3em]
\dot {\hat {\xrond }  } &=& f_{\xrond }(\hat y,\hat {\xrond })
\end{eqnarray*}

The remarkable fact here is that
there is no ``correction term'' in the dynamics of $\hat{\xrond}$.
Hence, we may expect that, if $P$ is a complete Riemannian metric for 
which there exist coordinates defined on some open set
$\Ouv  $
 satisfying (\ref{EQU9}), (\ref{EQU8}), and (\ref{LP57}) 
(with $\Ouv  $ replacing $\mathcal{N}_{x_0}$), then the system
\begin{equation}
\label{LP58}
\dot {\hat {\xrond }  } = f_{\xrond }(y,\hat {\xrond })
\end{equation}
(with $y$ instead of $\hat y$!) could be an appropriate
reduced order observer 
in charge of estimating the unmeasured components $\xrond$.
To show that this is indeed the case, we equip $\RR^{n-p}$, in which this 
reduced order observer lives, with 
the $y$ dependent Riemannian metric $\xrond\mapsto 
P_{\xrond\xrond}(y,\xrond)$. For each fixed $y$, we define the distance
\begin{equation}
\label{LP59}
d_\xrond(\xrond_a,\xrond_b;y)\!=\!\min_{\gamma_{_\xrond}}
\int_{s_a}^{s_b} \!\!\! \sqrt{
\frac{d\gamma _\xrond }{ds}(s)^\top
P_{\xrond\xrond}(y,\gamma _\xrond (s))\,  
\frac{d\gamma _\xrond }{ds}(s)}\,  ds
\end{equation}
where $\gamma_\xrond$ is any piecewise $C^1$ path satisfying
$
\gamma_\xrond (s_a)\;=\; \xrond_a$,
$\gamma_\xrond (s_b)\;=\; \xrond_b$.
With this, we have the following 
result for the reduced order observer (\ref{LP58}).

\begin{proposition}
\label{prop:ReducedGeodesicObserver}
Let $P_{\xrond\xrond}$ be a $y$-dependent Riemannian metric 
on $\RR^{n-p}$ and
$\mathcal{C}$  be
a closed subset of $\RR ^n$, with 
nonempty interior, 
satisfying
\begin{enumerate}
\item[A1'']:
$\mathcal{C}$ is weakly $P_{\xrond\xrond}$-geodesically convex in the following 
sense~: if $(\xrond_a,\xrond_b,y)$ is
such that $(y,\xrond_a)$ and $(y,\xrond_b)$ are in $\mathcal{C}$,
then there exists a minimizing geodesic
$[s_a,s_b]\ni s \mapsto \gamma _\xrond^*(s)$ in the sense of (\ref{LP59}) such that
$(y,\gamma _\xrond ^*(s))$ is in $\mathcal{C}$ for all $s$ in 
$[s_a,s_b]$.
Also, there exist coordinates denoted 
$\xrond$ and positive numbers
 $\plower , 
\pupper _{y1},
\bar f_{y1},$
such that, for each $(y,\xrond) $ in $ \mathcal{C}$, we have
$$
\begin{array}{c}
\displaystyle 
\plower  \,  I_{n-p}
\; \leq \; P_{\xrond\xrond}(y,\xrond)
\quad ,\qquad 
\left|\frac{\partial P_{\xrond\xrond}}{\partial y}(y,\xrond)\right|\; \leq \; 
\pupper _{y1} \\ \displaystyle
\left|\frac{\partial f_y}{\partial \xrond}(y,\xrond)\right|\; \leq \; \bar f_{y1}
\end{array}
$$
\item[A2'']:
There exists a strictly positive real number ${q}$ such that
(\ref{LP57}) holds on $\mathcal{C} \times  \SS^{n-p-1}$.
\end{enumerate}
Then, along the solutions to the system
$$
\dot y = f_y(y,\xrond )
\quad ,\qquad 
\dot {\xrond } = f_\xrond (y,\xrond )
\quad ,\qquad 
\dot {\hat {\xrond }  } = f_{\xrond }(y,\hat {\xrond })
\  ,
$$
the following holds:
$$
\Did d_\xrond ({\hat{\xrond }},\xrond;y  )
\; \leq \; 
-\underline{r}
\,  d_\xrond({\hat{\xrond }},\xrond ;y ) 
\  ,
$$
for all $(\xrond ,{\hat{\xrond }},y)$ such that
$(y,\xrond), (y,\hat{\xrond }) \in \mathcal{C}$ and
\begin{equation}
\label{LP60}
d_\xrond ({\hat{\xrond }},\xrond  )\; \leq \; 
\frac{(q-
2
\underline{r})\underline p\sqrt{\plower }}{ 
\bar p_{y1}
 \bar f_{y1}}
\  .
\end{equation}
\end{proposition}

The rationale is that, if the system is strongly 
differentially detectable
(see Definition~\ref{def:DifferentialDetectability}), 
then there exists a  reduced order observer 
that is exponentially convergent
as
long as
$(y,\xrond)$ and $(y,\hat \xrond)$ are in $\mathcal{C}$ and the
the coordinates 
$x=(y,\xrond)$ exist, which,
when $p=1$, 
we know is the case on a neighborhood 
of any point where $h$ has rank $1$.

\begin{proof}
Let $(\xrond ,\hat \xrond,y)$ be such that
$(y,\xrond )$ and $(y,\hat \xrond )$ are in $\mathcal{C}$.
From our assumption, there exists a minimizing geodesic
 $[s,\hat s] \ni s' \mapsto \gamma _\xrond^*(s')$ such that
$(y,\gamma _\xrond ^*(s'))$ is in $\mathcal{C}$ for all $s'$ in 
$[s,\hat s]$.
By following the same steps as in 
\cite[Proof of Theorem 2]{Aghanna.Rouchon.03} and with
\cite[(36)]{Sanfelice.Praly.11}, we can show that we have
\begin{eqnarray*}
\Did d_\xrond ({\hat{\xrond }},\xrond ;y )\; &\leq& 
\!\int_s^{\hat s}
\frac{
\!\!\frac{d\gamma _\xrond ^*}{ds}(r)^\top\!\!
\left[\mathcal{L}_{f_\xrond}P_{\xrond\xrond}(y,\gamma _\xrond^*(r))
+\frac{\partial P_{\xrond\xrond}}{\partial \xrond}(y,\gamma _\xrond^*(r))\,  
\dot y
\right]
\!\!\frac{d\gamma _\xrond ^*}{ds}(r)
}{
2\sqrt{
\frac{d\gamma _\xrond ^*}{ds}(r)^\top
P_{\xrond\xrond}(y,\gamma _\xrond^*(r))\,  
\frac{d\gamma _\xrond ^*}{ds}(r)
}
}\!  dr
\end{eqnarray*}
where
$
\dot y\;=\; f_y(y,\xrond)
\  .
$
So our result holds if the term between brackets is  upper bounded 
by $-2\underline{r} P(y,\gamma _\xrond^*(r))$.
Note that, in the coordinates given by A1'',
(\ref{LP57}) can be rewritten as
\begin{eqnarray}
\mathcal{L}_{f_\xrond} P_{\xrond\xrond}(y,\gamma _\xrond ^*)
\;+\; 
\frac{\partial P_{\xrond\xrond}}{\partial y}(y,\gamma _\xrond ^*)\,  \dot y
&\!\! \leq\!\!  &
-\qlower\,  P_{\xrond\xrond}(y,\gamma _\xrond ^*) 
\label{eqn:LfxrondP}
\;+\; 
\frac{\partial P_{\xrond\xrond}}{\partial y}(y,\gamma _\xrond ^*)
\left[f_y(y,\xrond)-f_y(y,\gamma _\xrond ^*)\right]
\end{eqnarray}
for all $(\xrond ,\gamma _\xrond ^*,y)$ such that $(y,\xrond)$ and 
$(y,\gamma _\xrond ^*)$ are in $\mathcal{C}$. But we have also
\begin{eqnarray*}\displaystyle 
\left|\frac{\partial P}{\partial y}(y,\gamma _\xrond ^*(r))
\left[f_y(y,\xrond)-f_y(y,\gamma _\xrond ^*(r)\right]\right|
\ & \leq & 
\  \displaystyle 
\bar p_{y1}
 \bar f_{y1}
\frac{d_\xrond ({\hat{\xrond }},\xrond ; y )}{\sqrt{\plower }}
\frac{P_{\xrond\xrond}(y,\gamma _\xrond ^*(r))}{\underline p}
\  .
\end{eqnarray*}
Hence, the result holds when \eqref{LP60} holds.
\end{proof}
\par\vspace{1em}\noindent
In this proof we see that the restriction (\ref{LP60}) disappears
and $q$ can be zero,
 if $\bar p_{y1}$ is zero, i.e.,~if $P_{\xrond\xrond}$ does not depend on $y$.
This is indeed the case when the level sets of the output function 
are totally geodesic as shown in \cite{Sanfelice.Praly.11}.
Hence, we have the following result.

\begin{proposition}
\label{LPprop1}
Under conditions A1'' and A2'' in Proposition~\ref{prop:ReducedGeodesicObserver}
with $q$ possibly zero,
if $P_{\xrond\xrond}$ does not depend 
on $y$, we have
\begin{equation}
\label{LP3}
\Did d_\xrond ({\hat{\xrond }},\xrond  )
\; \leq \; 
-q\,  d({\hat{\xrond }},\xrond  ) 
\end{equation}
for all $(\xrond ,{\hat{\xrond }},y) $ such that
$(y,\xrond)$ and $(y,\hat{\xrond })$ are in $\mathcal{C}$.
\end{proposition}

Again, the rationale is that if, the system is strongly
(respectively weakly)
differentially detectable and the output function level sets are totally 
geodesic, then there exists a reduced order observer which
makes the zero error set
$\{(y,\xrond,\hat \xrond):\,  \xrond=\hat  \xrond\}
$
exponentially stable (respectively stable) as
long as
$(y,\xrond)$ and $(y,\hat \xrond)$ are in $\mathcal{C}$ and the
coordinates 
$x=(y,\xrond)$ exist.

\begin{example}
\label{ex1}
Consider the harmonic oscillator with unknown frequency
(\ref{eqn:HarmonicOscillator}).
Its dynamics expressed in the coordinates $(y,\xrond_1,\xrond_2)$ 
we have obtained in (\ref{LP2}) are~:
\begin{equation}
\label{LP5}
\begin{array}{rcl}
\dot y & =&  \xrond_1 + k y, \\
\dot \xrond _1  & = & - y \,  (\xrond_2
-
\ell y^2)
\;-\; k\,  (\xrond_1 + k y), \\
\dot \xrond_2 & = & 
2 \ell y\,  (\xrond_1 + k y)
\end{array}
\end{equation}
In Example \ref{ex0}, we have shown this system is weakly 
differentially detectable with a metric the expression of which in 
the $(y,\xrond_1,\xrond_2)$ coordinates is
\begin{equation}
\label{LP4}
P((y,\xrond _1,\xrond _2))
\begin{array}[t]{@{}l@{}}
\begin{array}[t]{cl@{}}
=&\displaystyle 
\left[\left[\frac{\partial \varphi }{\partial x}(x)\right]^{-1}\right]^\top
P(x)
\left[\frac{\partial \varphi }{\partial x}(x)\right]^{-1}
\end{array}
\\
\begin{array}{cl}
=&
\left(\begin{array}{ccc}
1 
 & 0 & 
0
\\
0 & 2\ell & 0
\\
0
 & 0 & 1
\end{array}\right)
\end{array}
\end{array}
\end{equation}
As already observed in Example \ref{ex2},
 the decomposition given in (\ref{EQU8}) of
Theorem~\ref{thm:Eisenhart2} 
with even the 
$P_{\xrond\xrond}$ block independent of $y$.
So the assumptions of Proposition~\ref{LPprop1} are satisfied
with $\mathcal{C}=\RR^3$, but 
with $q=0$ and
the zero error set (with $\Ouv$ given in (\ref{LP7}))
$$
\mathcal{Z}\;=\; \{(y,\xrond _1,\xrond _2,\hat \xrond_1,\hat \xrond 
_2) \in \Ouv\times \RR^2:\: \xrond_2=\hat \xrond_2\}
$$
is globally stable.
To check that we have actually global stability, we note that the Lie derivative of the $P_{\xrond\xrond}$ 
block of $P$ in (\ref{LP4}) along the vector field given by 
(\ref{LP5})
satisfies for all $y$
$$
2\, \textsf{Sym}\left(
\vrule height 2em depth 2em width 0pt
\left(\begin{array}{cc}
2\ell & 0 \\ 0 & 1
\end{array}\right)
\left(\begin{array}{cc}
-k& -y \\ 
2
 \ell y & 0
\end{array}\right)
\right)
\;=\; 
\left(\begin{array}{cc}
-4\ell k & 
0 \\ 
0 & 0
\end{array}\right)
$$
where for a matrix $A$, $ \textsf{Sym}(A) = \frac{A+A^\top}{2}$.
This establishes that the vector field $f_\xrond$
defined as
$$
f_\xrond(y,\xrond)\;=\; 
\left(\begin{array}{@{}c}
- y \,  (\xrond_2
-
\ell y^2)
\;-\; k\,  (\xrond_1 + k y)
\\
2
 \ell y\,  (\xrond_1 + k y)
\end{array}
\right)
$$
is weakly geodesically monotonic uniformly in $y$. This 
implies that the flow it generates is a weak contraction. The solutions 
of the harmonic oscillator being bounded, the same holds for the 
solutions of
\begin{equation}
\label{LP83}
\dot {\hat {\xrond }}\;=\; f_\xrond(y,\hat \xrond)
\end{equation}
Then, according to \cite[Theorem 2]{Forni.Sepulchre.14},
the set\footnote{%
This means that the initial condition for $(x_1,x_2)$ 
is not the origin.
}
$$
\mathcal{Z}\: \mbox{\Large$\setminus$} \left(
\vrule height 0.6em depth 0.6em width 0pt
\varphi \left(
\vrule height 0.41em depth 0.41em width 0pt
\{(0,0)\}\times\RR_{>0}
\right)\: \times \RR^2\right)
\  ,
$$
with $\varphi$ defined in (\ref{LP2}),
is globally asymptotically stable for the interconnected system 
(\ref{eqn:HarmonicOscillator}), (\ref{LP83}).
\hfill $\triangle$
\end{example}

\section{Design of Riemannian Metric $P$ for Linearly Reconstructible Systems}
\label{sec:Reconstructibility}

We have seen in \cite[Theorem 2.9]{Sanfelice.Praly.11}
(see also \cite[Proposition 3.2]{Sanfelice.Praly.09.CDC.Observers}) 
that differential detectability implies that
each linear (time 
varying) system given by the first order approximation of 
(\ref{eqn:Plant1}) (assumed to be forward complete) along any 
of its solution is uniformly detectable. In \cite[Proposition 3.2]{Sanfelice.Praly.09.CDC.Observers}
it is also shown that, if this uniform linear detectability is 
strengthened into a uniform reconstructibility property (or, say, uniform 
infinitesimal
observability \cite[Section I.2.1]{Gauthier.Kupka.01.Book.Observers}), then a Riemannian metric 
exhibiting differential detectability does exist.
In this section, we recover
this last property through the solution of a Riccati equation and propose
a numerical method to compute the metric $P$.\footnote{
Some of the material in this section is in
\cite{Sanfelice.Praly.15.CDC}, which we reproduce here for the sake of completeness.
}

To do all this, we assume the existence of a
backward invariant open set 
$\Ouv $
for the system
\eqref{eqn:Plant1}.
This implies that, for each $x$ in $\Ouv $,
there exists a strictly positive real number $\sigma _x$, possibly infinite, such that
the corresponding solution to \eqref{eqn:Plant1},
$t\mapsto \XX(x,t)$, 
is defined with values in 
$\Ouv $
over $(-\infty ,\sigma _x)$.
For each such $x$,
the linearization of $f$ and $h$ evaluated along
$t \mapsto \XX(x,t)$
gives the functions 
$
t\mapsto A_{x}(t)  \;=\;   \frac{\partial f}{\partial \x}(\XX(x,t))$
and
$t\mapsto C_{x}(t) \;=\; \frac{\partial h}{\partial \x}(\XX(x,t))$,
which are defined on $(-\infty ,\sigma _x)$.
To these functions, we associate the following family of linear time-varying systems
with state $\xi $ in $ \RR^n$ and output $\eta $ in $ \RR^p$:
\begin{equation}
\label{LP6}
\dot \xi \;=\;  A_{x}(t)\,  \xi
\  ,\quad 
\eta \;=\; C_{x}(t)\,  \xi,
\end{equation}
which is parameterized by the initial condition $\x$ of the
chosen solution $t \mapsto \XX(x,t)$.
Below,  $\Phi_x$ denotes the
state transition matrix  for \eqref{LP6}.
It satisfies
$$ 
\frac{\partial \Phi_x}{\partial s}(t,s)\;=\; A_x(t) \Phi_x(t,s)
\quad ,\qquad
\Phi_x(s,s)\;=\; I
\  .
$$

\begin{definition}[reconstructibility]
The family of systems \eqref{LP6} is said to be reconstructible
on a set 
$\Ouv $
if there exist strictly positive real numbers $\tau$ and $\varepsilon$ 
such that we have
\begin{equation}
\label{LP30}
\int_{-\tau } ^0
\Phi _{x}(t,0)^\top C_{x}(t)^\top C_{x}(t) \Phi _{x}(t,0) dt\; \geq \; \varepsilon  \,  I \qquad
\forall x \in \Ouv \  .
\end{equation}
\end{definition}

\begin{proposition}
\label{propo:NecessitySatisfactionReconstructibilityDRE}
Let $Q$ be a symmetric contravariant $2$-tensor.
Assume there exist
\begin{list}{}{%
\parskip 0pt plus 0pt minus 0pt%
\topsep 0pt plus 0pt minus 0pt
\parsep 0pt plus 0pt minus 0pt%
\partopsep 0pt plus 0pt minus 0pt%
\itemsep 0pt plus 0pt minus 0pt
\settowidth{\labelwidth}{ii)}%
\setlength{\labelsep}{0.5em}%
\setlength{\leftmargin}{\labelwidth}%
\addtolength{\leftmargin}{\labelsep}%
}
\item[i)]
an open set $\Ouv  \subset \RR^n$ that is
backward invariant 
for \eqref{eqn:Plant1}
and 
on which the family of systems \eqref{LP6} is 
 reconstructible;
\item[ii)]
coordinates for $x$ such that the derivatives of $f$ and $h$ are bounded on
$\Ouv $
and we have
\\[0.3em]\null \hfill $\displaystyle 
0 \; < \;  \underline{q}\,  I \; \leq \;  Q(\x)  \leq \;  \overline{q} \,  I
\qquad \forall x\in \Ouv \ .
$\hfill \null 
\refstepcounter{equation}\label{eqn:QboundsOnC}$(\theequation)$
\end{list}
Then, there exists a symmetric covariant 
$2$-tensor $P$
defined
on $\Ouv $,
which admits a Lie derivative $\mathcal{L}_f P$ satisfying
\begin{equation}
\label{eqn:PidotNL}
\mathcal{L}_f P(x) \;= \; 
d \h(\x) \otimes d \h(\x)  \;-\;   P(x) Q(x) P(x) \qquad \forall x \in \Ouv \ ,
\end{equation}
and there exist strictly positive real numbers $\plower$ and $\pupper$
such that, in the coordinates given above, we have
\begin{eqnarray}\label{eqn:PboundsOnC}
0 \; < \;  \plower\,  I \; \leq \;  P(\x)  \leq \;  \pupper \,  I
\qquad \forall x\in \Ouv \ .
\end{eqnarray}
\end{proposition}

\begin{proof}
The proof of Proposition~\ref{propo:NecessitySatisfactionReconstructibilityDRE}
can be found in \cite{Sanfelice.Praly.15.CDC}. It
relies on a fixed point argument, the core of which is the fact the
flow generated by the differential Riccati equation is a contraction. This
fact, first established for the discrete time case in \cite{Bougerol.93.SIAM},
is proved in \cite{bonnabel2013geometry} for the continuous-time
case. 
\end{proof}

\begin{remark}~
In his introduction of Riccati differential equations 
for matrices in \cite{Rad1927,Rad1928}, Radon has shown that such equations can be solved via two coupled linear 
differential equations.  (See also \cite{Reid72}.) In our framework, this leads to 
obtain a 
solution to equation (\ref{eqn:PidotNL}) by solving in $(\alpha,\beta)$ 
the coupled system
\begin{equation}
\label{LP90}
\begin{array}{rcl}
\displaystyle \sum_{i=1}^n \frac{\partial \alpha}{\partial x_i}(x) f_i(x)
&=&\displaystyle -\frac{\partial f}{\partial x}(x)^\top \alpha (x) 
+ \displaystyle
\frac{\partial h}{\partial x}(x) \frac{\partial h}{\partial x}(x)^\top 
\beta(x)
\  ,
\\[0.5em]
\displaystyle \sum_{i=1}^n \frac{\partial \beta }{\partial x_i}(x) f_i(x)&=&
\displaystyle 
Q(x) \alpha (x) + \frac{\partial f}{\partial x}(x) \beta(x)
\end{array}
\end{equation}
with $\beta$ invertible
and then picking
$
P(x)\;=\; \alpha (x) \,  \beta (x)^{-1}
$.
\hfill $\Box$
\end{remark}

\begin{remark}
\label{remark:ComparisonToKalman}
Our observer 
in \eqref{LP56}
with right-hand side given by
\eqref{eqn:GeodesicObserverVectorField} or \eqref{4} resembles the Extended Kalman filter
for a particular choice of $\delta$.
In fact,
when the metric is obtained by solving (\ref{eqn:PidotNL}), the observer we obtain from 
(\ref{eqn:GeodesicObserverVectorField}) (or \eqref{4})
with $\delta (y_a,y_b)\;=\; |y_a-y_b|^2$ resembles an Extended Kalman Filter (see 
\cite{jazwinski1970stochastic} for instance) since, in
some coordinates, our observer is
\begin{eqnarray}
\label{LP84}
\dot \xhat &\!\!\hskip -0.3em =&\!\!\hskip -0.3em 
 f(\xhat) \!\!\;-\;\!\! 2\, k_E(\xhat)\,  P(\xhat)^{-1} \frac{\partial h}{\partial x}(\xhat)^\top
\left(\h(\xhat) - y\right),\ \ 
\end{eqnarray}
\begin{eqnarray}
\sum_{i=1}^n \frac{\partial P}{\partial x_i}(\hat x) f(\xhat)&\hskip 
-0.3em\!\! =&\hskip -0.3em\!\! 
- P(\hat x) \frac{\partial f}{\partial x}(\hat x)
-
\frac{\partial f}{\partial x}(\hat x)^\top P(\hat x) 
\label{LP85}
+\frac{\partial h}{\partial x}(\hat x)^\top\frac{\partial h}{\partial x}(\hat x)
- P(\hat x) Q(\hat x) P(\hat x)
\qquad \null 
\end{eqnarray}
while the corresponding extended Kalman filter would be
\begin{eqnarray}
\label{LP86}
\dot \xhat &=& f(\xhat) \;-\;   P^{-1} \frac{\partial h}{\partial x}(\xhat)^\top
\left(\h(\xhat) - y\right)
\  ,
\\
\dot P&=&
- P \frac{\partial f}{\partial x}(\hat x)
-
\frac{\partial f}{\partial x}(\hat x)^\top P 
\label{LP87}
+\frac{\partial h}{\partial x}(\hat x)^\top\frac{\partial h}{\partial x}(\hat x)
- P Q(\hat x) P\  .
\end{eqnarray}
The expressions for $\dot{\hat x}$ in (\ref{LP84}) and (\ref{LP86}) are the same except for the 
presence of $k_E$ in (\ref{LP84}). On the other hand, (\ref{LP85}) and 
(\ref{LP87}) are significantly different. The former is a 
partial differential equation which can be solved off-line as
an algebraic Riccati equation. If the assumptions in
Proposition~\ref{propo:NecessitySatisfactionReconstructibilityDRE} 
are satisfied, (\ref{LP85}) has a solution, guaranteed to be bounded and 
positive definite
on $\Ouv$.
Nevertheless, assumption 
A3 of Theorem \ref{thm1part1} may not hold
but then according to Proposition \ref{prop1}, we have
a locally convergent observer.

The differential Riccati 
equation (\ref{LP87}) of the extended Kalman filter is an ordinary 
differential equation with $P$ being part of the observer state.
The corresponding observer is also known to be locally convergent
but under the extra assumption that $P$ is bounded and 
positive definite. See \cite{Bonnabel-Slotine.15} for instance. 
Unfortunately, 
even when the assumptions in
Proposition~\ref{propo:NecessitySatisfactionReconstructibilityDRE} 
are satisfied, we have no guarantee that $P$ has such properties
except may be if $\hat x$ remains close enough to 
$x$ (which is what is to be proved).
\hfill $\Box$
\end{remark}

The quadratic term
$P(x) Q(x) P(x)$ in
the ``algebraic Riccati equation'' \eqref{eqn:PidotNL}, can be 
replaced by 
$\lambda P(x)$. Specifically we have the following reformulation of 
\cite[Proposition 3.2]{Sanfelice.Praly.09.CDC.Observers}.

\begin{proposition}
\label{prop:PidotLambda}
Under the conditions of Proposition~\ref{propo:NecessitySatisfactionReconstructibilityDRE},
there exists
$\underline{\lambda} > 0$ such that, 
for each $\lambda > \underline{\lambda}$,
there exists
a symmetric covariant 
$2$-tensor $P$
defined
on $\Ouv $
that admits a Lie derivative $\mathcal{L}_f P$ satisfying
\begin{equation}
\label{eqn:PidotNLlambda}
\mathcal{L}_f P(x) \;= \; 
d \h(\x) \otimes d \h(\x)  \;-\;   \lambda P(x)  \qquad \forall x \in \Ouv \ ,
\end{equation}
and
there exist
strictly positive real numbers $\plower$ and $\pupper$
such that the expression of $P$ in the coordinates given by the 
assumption satisfies \eqref{eqn:PboundsOnC}.
\end{proposition}

\begin{proof}
See \cite{Sanfelice.Praly.15.CDC}.
\end{proof}

\begin{remark}
When the metric is given by (\ref{eqn:PidotNLlambda}), the observer 
we obtain from
(\ref{eqn:GeodesicObserverVectorField})  with $\delta (y_a,y_b)\;=\; |y_a-y_b|^2$ 
resembles the Kleinman's observer, dual of the 
Kleinman's controller proposed in \cite{Kleinman.70}.
Indeed, in 
some coordinates, our observer is
\begin{eqnarray*}
\dot \xhat \!\!\!&=&\!\!\!
 f(\xhat) \;-\; 2\, k_E(\xhat)\,  
P(\hat x) ^{-1}
\frac{\partial h}{\partial x}(\xhat)^\top
\left(\h(\xhat) - y\right)
\  ,
\\
P(x)\!\!\! & =&\!\!\!
\lim_{T\to \infty }\int_{-T}^0 \!\!
\exp(\lambda t)
\Phi _{x} (t,0)^\top C_{x}(t)^\top C_{x}(t) \Phi _{x} (t,0) 
dt,
\end{eqnarray*}
the latter being a solution to (\ref{eqn:PidotNLlambda}). 
Correspondingly,
 Kleinman's observer would be
\begin{eqnarray*}
\dot \xhat &=& f(\xhat) \;-\;   W(\xhat)^{-1} \frac{\partial h}{\partial x}(\xhat)^\top
\left(\h(\xhat) - y\right)
\  ,
\\
W(x)&=&
\int_{-T}^0
\Phi _{x} (t,0)^\top C_{x}(t)^\top C_{x}(t) \Phi _{x} (t,0) 
\,  dt
\end{eqnarray*}
with $T$ positive.
\hfill $\Box$
\end{remark}

\begin{example}
\label{ex:RevisitHarmonicOscillator-ComputationofP}
For the harmonic oscillator with unknown 
frequency (\ref{eqn:HarmonicOscillator}), it can be checked that
the following expression of $P$ 
is a solution to (\ref{eqn:PidotNLlambda}):
\\[1em]
$\displaystyle 
P(x)\;=\; \left(\begin{array}{c@{\  ,\  }c@{\  ,\  }c}
\displaystyle 
\frac{\lambda ^2 +2x_3 }{\lambda (\lambda ^2+4x_3)}
&\displaystyle 
\star
&
\displaystyle 
\star\\[1.56em]
\displaystyle 
-\frac{1  }{ (\lambda ^2+4x_3)}
&\displaystyle 
\frac{2}{\lambda (\lambda ^2+4x_3)}
&
\displaystyle 
\star\\[1.56em]
\displaystyle 
\frac{-
\lambda ^3
x_1
+
(\lambda ^2-4x_3)
x_2}{\lambda ^2(\lambda ^2+4x_3)^2}
&
\displaystyle 
\frac{ (3\lambda ^2+4x_3)
x_1
-
4\lambda 
x_2
}{\lambda ^2(\lambda ^2+4x_3)^2}
&
a
\end{array}\right)
$\refstepcounter{equation}\label{LP97}\hfill$(\theequation)
$\\[1em]
where the various $\star$ should be replaced by their symmetric 
values and
\begin{eqnarray*}
a &=&
\displaystyle 
\frac{6\lambda ^4
+12\lambda ^2x_3
+16x_3^2
}{\lambda ^3(\lambda ^2+4x_3)^3}x_1 ^2
-
\frac{
4(5\lambda ^2+4x_3)
}{\lambda ^2(\lambda ^2+4x_3)^3}
x_1x_2
+
\frac{
4(5\lambda ^2+4x_3)
}{\lambda ^3(\lambda ^2+4x_3)^3}
x_2^2 
\end{eqnarray*}
\hfill $\triangle$
\end{example}
\par\vspace{1em}
One
way to prove Proposition
\ref{propo:NecessitySatisfactionReconstructibilityDRE}, 
respectively Proposition \ref{prop:PidotLambda}, is to show that
the system
$$
\begin{array}{rcl}
\dot{x} & = & f(x)
\  ,
\\[0.5em]
\dot{\pi} &=&\displaystyle 
F(x,\pi )\;=\; 
-\pi \,
\frac{\partial f}{\partial x}(x)
 - 
\frac{\partial f}{\partial x}(x)^\top\, \pi 
\displaystyle + 
\frac{\partial h}{\partial x}(x)^\top 
\frac{\partial h}{\partial x}(x)
- \pi\, Q(x)\, \pi 
\  ,
\end{array}
$$
respectively
$$
\begin{array}{rcl}
\dot{x} & = & f(x)
\  ,
\\[0.5em]
\dot{\pi} &=&\displaystyle 
F(x,\pi )\;=\; 
-\pi \,
\frac{\partial f}{\partial x}(x)
 - 
\frac{\partial f}{\partial x}(x)^\top\, \pi 
\displaystyle + 
\frac{\partial h}{\partial x}(x)^\top 
\frac{\partial h}{\partial x}(x)
- 
\lambda \,  \pi
\  ,
\end{array}
$$
admits an invariant manifold of the form
$\{(x,\pi)\, : \,  \pi=P(x)\}$.
These facts suggest the following method to approximate $P$.\\
Given $x$ in $\Ouv$ at which $P$ is to be evaluated, 
pick $T > 0$ large enough, 
and 
perform the following steps\footnote{%
In the case where the system is time varying and its time variations 
are dealt with as explained in footnote \ref{footnote1}, these steps 
do require the knowledge of the time functions. This imposes a  
difficulty when,  for instance, the time functions are induced by inputs provided by a feedback 
law.
}:
\begin{list}{}{}
\item[\bf Step 1)] Compute the solution $[-T,0] \ni t \mapsto \XX(x,t)$ to \eqref{eqn:Plant1} 
backward in time from the initial condition $\x$ at time $t=0$, 
up to a negative time $t=-T$;
\item[\bf Step 2)]
Using the function $[-T,0] \ni t \mapsto \XX(x,t)$ obtained 
in Step 1,
compute the solution $[-T,0] \ni t 
\mapsto \Pi(t)
$
with initial condition
$\pi(-T) = \plower \,  I_n$,
to 
\begin{eqnarray*}
\dot \pi  &=&  
-\pi \frac{\partial f}{\partial \x}(\XX(x,t)) \;-\; \frac{\partial f}{\partial \x}(\XX(x,t))^\top \pi
\;+\; \frac{\partial h}{\partial \x}(\XX(x,t))^\top \frac{\partial h}{\partial \x}(\XX(x,t))
\;-\; \pi Q(X(x,t)) \pi
\  ,
\end{eqnarray*}
respectively to
\begin{eqnarray*}
\dot \pi  & = &
-\pi \frac{\partial f}{\partial \x}(\XX(x,t)) \;-\; \frac{\partial f}{\partial \x}(\XX(x,t))^\top \pi 
\;+\; \frac{\partial h}{\partial \x}(\XX(x,t))^\top \frac{\partial h}{\partial \x}(\XX(x,t))
\;-\; \lambda \,  \pi 
\end{eqnarray*}
with $\lambda $ large enough.
\item[\bf Step 3)] Define the value of $P$ at $x$ as the value $\Pi(0)$.
\end{list}
\medskip

By griding the state space of $x$ and approximating $P$ at each such $x$,
the method suggested above can be considered as a design tool, at least for low 
dimensional systems.  
Note that the computations in Step 1 and Step 2 only require 
the use of a scheme for integration of ordinary differential 
equations. 
In the following example, we employ this method to approximate the
metric $P$ for the harmonic oscillator after a convenient reparameterization 
allowing a reduction of
 the number of points needed
in a grid for a given desired precision.

\begin{example}
\label{ex:RevisitHarmonicOscillator-NumericalComputationofP}
The second version of the
proposed algorithm applied to the harmonic oscillator in \eqref{eqn:HarmonicOscillator}  
leads to an approximation of the analytic expression of the metric $P$ given 
in Example~\ref{ex:RevisitHarmonicOscillator-ComputationofP}.
To this end, we exploit
the fact that
$\sqrt{x_3}$ and $t$ have the same dimension
and, similarly,
$\frac{x_1}{\sqrt{x_3}}$, and $\frac{x_2}{x_3}$  have the same 
dimension.
To exploit this property, we let
$$
r\;=\; \sqrt{x_3x_1^2+x_2^2}
\quad ,\qquad 
\cos(\theta )\;=\; \frac{\sqrt{x_3}x_1}{r}
$$
$$\sin(\theta )\;=\; \frac{x_2}{r}
\quad ,\qquad 
\omega \;=\; \frac{\lambda }{\sqrt{x_3}}
\  .
$$
Then, it can be checked that the metric $P$ can be factorized as
$$P(x_1,x_2,x_3,\lambda ) =
M(x_3)^{-1}P\left(\cos(\theta ),\sin(\theta ),1,\,  \omega\right)
M(x_3)^{-1}\ ,$$
where
$
M(x_3) = \mbox{diag}\left(
x_3^{1/4}, \sqrt{x_3}x_3^{1/4}, \frac{x_3\textstyle\sqrt{x_3}x_3^{1/4}}{r}\right)
$.
This shows that 
it is sufficient to know the function
$(\theta ,\omega )\ni (\SS^1\times\RR_{>0})\mapsto
P(\cos(\theta ),\sin(\theta ),1,\omega )$ 
 and the value of $x_3$
to know the function $P$ everywhere
on $(\RR^2\setminus\{0\})\times \RR_{>0}^2$.
Further using the fact that
$$
\frac{\partial h}{\partial x}(x_1,x_2,x_3)\;=\; 
\left(\begin{array}{ccc}
1 & 0 & 0
\end{array}\right)
\  ,
$$
the gain of the proposed observer reduces to
\\[1em]\vbox{\noindent
$\displaystyle 
P(x_1,x_2,x_3,\lambda )^{-1}\frac{\partial h}{\partial x}(x_1,x_2,x_3)^\top
=
\left(\begin{array}{ccc}
\displaystyle 
\sqrt{x_3}& 0 & 0
\\
0 &  
\displaystyle x_3  & 0
\\
0 & 0 &\displaystyle \frac{x_3^2}{r}
\end{array}\right)
P\left(\cos(\theta ),\sin(\theta ),1,\,  \omega\right)^{-1}
\left(\begin{array}{@{\,  }c@{\,  }}
1 \\[0.6em] 0 \\[0.6em] 0
\end{array}\right)
$}\\[1em]
This shows  that 
it is sufficient to know the function
$$
(\theta ,\omega )\ni (\SS^1\times\RR_{>0})
\; \mapsto\; 
\left(
\begin{array}{c}
P^{-1}_{11}(\cos(\theta ),\sin(\theta ),1,\omega )
\\
P^{-1}_{12}(\cos(\theta ),\sin(\theta ),1,\omega )
\\
P^{-1}_{13}(\cos(\theta ),\sin(\theta ),1,\omega )
\end{array}\right)
$$
to know the observer gain everywhere
on $(\RR^2\setminus\{0\})\times \RR_{>0}^2$.
Hence it is sufficient to grid the circle $\SS^1$ with $m_\theta $ points and the 
strictly positive real numbers with $m_\omega $ points, 
and therefore to
store only $3*m_\theta *m_\omega $ values in which the above function is 
interpolated.

We note that for the  computation of $P$ using the algorithm above, 
since a closed-form expression of the solutions to \eqref{eqn:HarmonicOscillator} is available, Step 1 of the algorithm is not needed.
To compute an approximation of $P$, 
we define a grid of the $(\theta,\omega)$-region $[-\pi,\pi]\times[4,7]$ with  $m_{\theta} * m_{\omega}$ points
with $m_{\theta}=360$ and $m_{\omega}=100$.
The value of $T$ used in the simulations is chosen as a function of $\omega$, namely, $T(\omega)$,
so as to guarantee a desired absolute error for the approximation of $P$
for the given point $(\theta,\omega)$ from the grid.

\begin{figure}[h!] 
\begin{center}
\psfrag{2.8}[][][0.9]{}
\psfrag{2.6}[][][0.9]{}
\psfrag{2.4}[][][0.9]{}
\psfrag{2.2}[][][0.9]{}
\psfrag{1.8}[][][0.9]{}
\psfrag{1.6}[][][0.9]{}
\psfrag{1.4}[][][0.9]{}
\psfrag{1.2}[][][0.9]{}

\psfrag{x1}[][][0.9]{$x_1, \hat{x}_1$}
\psfrag{x2}[][][0.9]{$x_2, \hat{x}_2$}
\psfrag{x3}[][][0.9]{$x_3, \hat{x}_3$}
\psfrag{e1}[][][0.9]{$x_1 - \hat{x}_1$}
\psfrag{e2}[][][0.9]{$x_2 - \hat{x}_2$}
\psfrag{e3}[][][0.9]{$x_3 - \hat{x}_3$}
\psfrag{t}[][][0.6]{$t [s]$}

\subfigure[Solutions. \label{fig:3D}]
{\includegraphics[trim=2cm 1cm 1cm 1cm,clip,width=.52\textwidth]{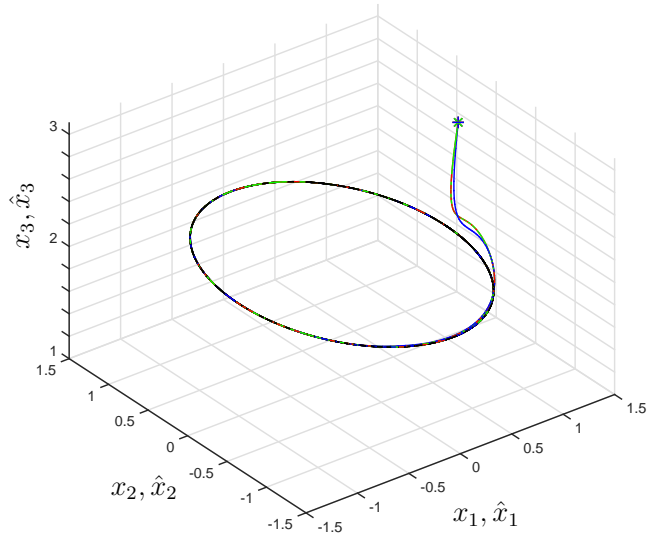}} \qquad\qquad
\subfigure[Estimation errors. \label{fig:Errors}]
{\includegraphics[trim=1.3cm 0.2cm 1cm 0.2cm,clip,width=.49\textwidth]{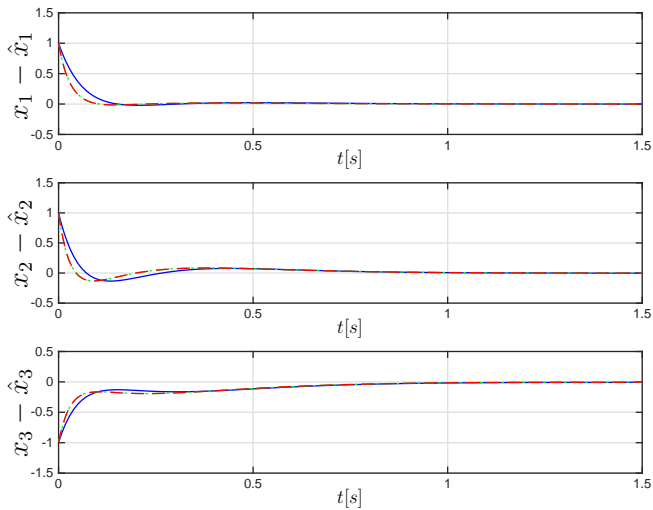}} 
\caption{Solutions to the observer converging to the estimate obtained with exact gain with $\lambda=8$ (solid blue/darkest), with exact gain discretized over a grid (dash dot blue/gray), and with computed and interpolated gain (dashed red/dark).}
\label{fig:ObserverSolutionsConverging}
\end{center}
\end{figure}

Figure~\ref{fig:ObserverSolutionsConverging} shows 
state estimates $\hat{x}$ 
using the observer in \eqref{4} for a periodic 
solution to \eqref{eqn:HarmonicOscillator}.
These solutions start from the same initial condition 
and are such that the state estimates asymptotically converge to the periodic solution.
The solid blue/darkest solution corresponds 
to the estimate obtained using in \eqref{4} the analytic expression of
$P$ in \eqref{LP97} with parameter $\lambda = 8$, which 
is a large enough value to satisfy the desired precision.
The other solutions in Figure~\ref{fig:ObserverSolutionsConverging} correspond to estimates obtained with different computed values
of $P$ using our algorithm.  
The dash dot blue/gray solution is obtained when 
observer gain is discretized over the 
chosen grid and provided to the observer using nearest point interpolation.
The dashed red/dark solution is obtained when the observer gain is computed (over the same grid) using the algorithm proposed above.
For each simulation, the error trajectories converge to zero.
Note that the error between the
dash dot blue/gray solution and the
dashed red/dark solution is quite small.
As the figures suggest,
the estimates obtained with the approximated gains are close to
the one obtained with its analytical expression. Additional numerical analysis confirms that the error 
between the solutions gets smaller as the
number of points and the quality of the interpolation are increased.
\hfill $\triangle$
\end{example}

\section{Design of Riemannian Metric $P$ for Strongly 
Differentially Observable Systems}
\label{sec:High gain}

According to \cite[Definition 4.2 of Chapter 2]{Gauthier.Kupka.01.Book.Observers}, 
the nonlinear system \eqref{eqn:Plant1}
is strongly 
differentially observable of order $n_o$ on an open set $\Ouv $ if,
for the positive integer $n_o$,
the function $\euH_{n_o}:\Ouv  \to \RR^{m\times n_o}$ defined  as
\begin{equation}\label{eqn:StronglyDiffObservableMap}
\euH_{n_o}(x)\;=\; \left(\begin{array}{c}
h(x), 
L_f
h(x), \cdots, 
L_f
^{n_o-1}h(x)
\end{array}\right)^\top
\end{equation}
is an injective immersion, i.e., an injective map whose differential 
is injective at each point $x $ in $ \Ouv$.

\begin{example}
For the system (\ref{eqn:HarmonicOscillator})
in Example~\ref{ex:RevisitHarmonicOscillator-ComputationofP}, 
successive derivatives of $y$ lead to
$$
\euH_{3}(x)\;=\; \left(\begin{array}{@{}c@{}}
x_1, x_2, -x_3\,  x_1
\end{array}\right)^\top$$
$$
\euH_{4}(x)\;=\; \left(\begin{array}{@{}c@{}}
x_1, x_2, -x_3\,  x_1, \ - x_3 x_2
\end{array}\right)^\top
\  .
$$
The map $\euH_3$ is an injective immersion on $\Ouv _3 = 
\left(\RR\setminus\{0\}\right)\times\RR\times\RR_{>0}$ which is not an 
invariant set. Instead 
$\euH_4$ is an injective immersion on $\Ouv 
_4=\left(\RR^2\times\RR_{>0}\right)\setminus\left(\{(0,0)\}\times\RR_{+}\right)$ which is an invariant set.
Hence, the system in Example~\ref{ex:RevisitHarmonicOscillator-ComputationofP} is strongly 
differentially observable of order $4$ on the invariant set $\Ouv _4$.
\hfill $\triangle$
\end{example}

The property that $\euH_{n_o}$ is an injective immersion implies that 
the family of systems \eqref{LP6} is reconstructible (on $\Ouv$).
According to Section \ref{sec:Reconstructibility}, 
this property further implies that differential detectability holds with
a metric obtained as a solution of (\ref{eqn:PidotNL}) or of
(\ref{eqn:PidotNLlambda}).
But we can take advantage of the strong observability property to 
give another more explicit expression for the metric.
Precisely,
 we assume the following properties.
\begin{enumerate}
\item[B~:] 
There are coordinates for $x$ in $\Ouv$ such 
that
\begin{itemize}
\item
$\euH_{n_o}$ is Lipschitz and a uniform immersion, i.e., assume the existence of strictly positive real numbers 
$\underline{h}$ and $\overline{h}$ such that we have
\begin{equation}\label{eqn:euHassumption}
\underline{h}\,  I\; \leq \; \frac{\partial \euH_{n_o}}{\partial x}(x)^\top  \frac{\partial \euH_{n_o}}{\partial x}(x)
\; \leq \; \overline{h} \,  I
\qquad \forall x \in \Ouv
\  ;
\end{equation}
\item
There exists a strictly positive real number $\nu$
such that, in the given coordinates for $x$, we have the following Lipschitz-like 
condition\footnote{%
We say that (\ref{13}) is a Lipschitz-like condition since, 
the function $\euH_{n_o}$, being injective, has a left inverse 
$\euH_{n_o}^{\mbox{\tiny li}}$ satisfying
$
\euH_{n_o}^{\mbox{\tiny li}}\left(\euH_{n_o}(x)\right)\;=\; x$.
Consequently, we have
$
L _f^{n_o}h\;=\; L _f^{n_o}(h\circ\euH_{n_o}^{\mbox{\tiny li}}\circ\euH_{n_o})
$.
It follows that, if the function
$\ox\mapsto L _f^{n_o}(h\circ\euH_{n_o}^{\mbox{\tiny li}})(\ox)$
is Lipschitz, then (\ref{13}) 
holds.
}
\begin{equation}
\label{13}
\left|\frac{\partial L _f^{n_o}h}{\partial x}(x)\right|
\; \leq \; 
\frac{1}{\nu}
\,  \left|\frac{\partial \euH_{n_o}}{\partial x}(x)\right|
\  .
\end{equation}
\end{itemize}
\end{enumerate}
To exploit these properties, we note first that we have
\begin{eqnarray}
\nonumber
L_f\euH_{n_o}(x)
\!\!\!& = &\!\!\!
\frac{\partial \euH_{n_o}}{\partial x}(x)\,  f(x)\;=\; \euA\,  
\euH_{n_o}(x)\;+\; \euB \, 
L_f
^{n_o}h(x),\\ \label{LP61}
y & = & h(x)\;=\; \euC\,  \euH_{n_o}(x)
\end{eqnarray}
where $\euA$, $\euB$, and $\euC$ are given by
$$
\euA\;=\; \left(\begin{array}{ccccc}
0          & I_m          & 0           & \dots & 0
\\
\vdots & \ddots & \ddots & \ddots& \vdots
\\
\vdots & & \ddots & \ddots & 0
\\
\vdots & &  & \ddots & I_m
\\
0 & \dots & \dots & \dots & 0
\end{array}\right)
\quad ,\qquad 
\euB\;=\; \left(\begin{array}{c}
0 \\ \vdots \\ \vdots \\ 0 \\ I_m
\end{array}\right)$$
$$\euC\;=\; \left(\begin{array}{ccccc}
I_m & 0 & \dots & \dots & 0
\end{array}\right)
\  .
$$
Then,
among the many results known about high gain observers, 
we have the following property.
\normalcolor
\begin{lemma}
\label{lemma:BoundedRealLemma}
Given $\nu$ satisfying (\ref{13}), there exist an $(m\times n_o)\times(m\times n_o)$
symmetric positive definite matrix $\euPL $, a $(m\times n_o)
\times m$ column vector $\euKL $, 
and a strictly positive real number $q$ satisfying
\begin{equation}
\label{12}
\begin{array}{lll}
\euPL \left(\euA -\euKL \euC\right)
\;+\; 
\left(\euA -\euKL \euC\right)^\top \euPL 
\displaystyle \;+\; 2q\,  I_{m\times n_o}
\;+\; \frac{1}{q \nu^2} \euPL \,   \euB\,  \euB ^\top 
\euPL \; \leq \; 0
\  .
\end{array}
\end{equation}
\end{lemma}

With Lemma~\ref{lemma:BoundedRealLemma}, we pick $P$ as the metric induced by the 
immersion $\euH_{n_o}$. (See \cite[Example 2 of Chapter II]{Sakai.96}.)
Namely, in the coordinates $x$ given by assumption B 
so that (\ref{eqn:euHassumption}) and (\ref{13}) hold,
we express $P$ on $\Ouv$ as
\begin{equation}
\label{14}
P(x)\;=\; \frac{\partial \euH_{n_o}}{\partial x}(x)^\top \euPL\,  \frac{\partial \euH_{n_o}}{\partial x}(x)
\  .
\end{equation}
\begin{remark}
The above design of $P$ relies strongly on the high gain observer 
technique. Nevertheless, the observer 
we obtain differs from a usual high gain observer, at least when
$n_o$ is strictly larger than $n$, i.e.,
$\euH_{n_o}$ is an injective immersion and not a diffeomorphism. Indeed, the state $\hat x$ 
of our observer lives in $\RR^n$, whereas the state of a usual high 
gain observer would live in $\RR^{n_o}$, not diffeomorphic to $\RR^n$, and 
a left inverse of $\euH_{n_o}$ would be needed to extract $\hat x$ from this 
state.
\end{remark}

\begin{proposition}
\label{prop:LfPforSystemsWithStrongDiffObs}
Suppose that, with $\euH_{n_o}$ defined in
(\ref{eqn:StronglyDiffObservableMap}), Assumption B holds and let  
$\euPL $ be any symmetric positive definite matrix satisfying 
(\ref{12}). Then, \eqref{14} defines a positive definite symmetric 
covariant $2$-tensor which 
satisfies the differential detectability property (\ref{eqn:LfPtangentially}) on
$\Ouv$. 
\end{proposition}

Here, similar to Ackerman's formula for linear systems, where the observer gain
uses the inverse of the observability matrix, the gain of our observer, namely,  
 $P(\hat x)^{-1}\frac{\partial h}{\partial x}(\hat x)^\top$, resulting from expressing the metric as in (\ref{14})
is obtained by writing the system in an observable form.
This form can be obtained  using $\frac{\partial \euH_{n_o}}{\partial x}(x)$
as the observability matrix,
the inverse of which
also appears in the gain of our observer.

\begin{proof} We proceed by establishing the needed properties for $P$.

\noindent
$\bullet$\ \textit{$P$ is a symmetric covariant $2$-tensor~:}
Let $\widetilde x$ be other coordinates related to $x$ by
$
\widetilde x\;=\; \varphi(x)
$
with $\varphi$ being a diffeomorphism.
Let also $\widetilde h$, $\widetilde P$, and ${\widetilde 
\euH}_{n_o}$ denote the expression of $h$, $P$, and ${\euH}_{n_o}$ in the coordinates $\widetilde x$, respectively.
They satisfy
\begin{eqnarray*}
\widetilde h(\widetilde x)=h(x)
\  ,
\qquad
\widetilde f(\widetilde x)\ = \ \frac{\partial \varphi}{\partial x}(x)\,  f(x)
\ \qquad \\
 \frac{\partial h}{\partial x}(x)\ =\ \frac{\partial \widetilde h}{\partial 
\widetilde x}(\widetilde x)\frac{\partial \varphi}{\partial x}(x)
\  , \qquad
\euH_{n_o}(x)={\widetilde \euH}_{n_o}(\widetilde x)\,  \frac{\partial \varphi}{\partial x}(x)
\  \\
P(x)\ =\ \frac{\partial \varphi}{\partial x}(x)^\top \widetilde P(\widetilde x)\frac{\partial \varphi}{\partial x}(x)
\   \qquad\qquad\qquad
\end{eqnarray*}
the latter showing that $P$ satisfies the rule a linear operator should obey under a change of 
coordinates to be a symmetric covariant $2$-tensor.
\par\vspace{1em}\noindent
$\bullet$\ \textit{$P$ is positive definite~:}
Using \eqref{14} and the positive definiteness of $\euPL$, we have
$$
0\; <\; \lambda _{\min}(\euPL)\,  \underline{h}\,  I\; \leq \; P(x)\; \leq \; 
\lambda _{\max}(\euPL)\,  \overline{h}\,  I
\qquad \forall x \in \Ouv  
\  .
$$
\par\vspace{1em}\noindent
$\bullet$\ \textit{$P$ satisfies (\ref{eqn:LfPtangentially})~:}
With (\ref{LP61}) and (\ref{LP64}), we obtain
\\[1em]$\displaystyle 
\mathcal{L} _fP(x)
\;=\;
\frac{\partial \euH_{n_o}}{\partial x}(x)^\top
\left(
\euPL \,  \euA + \euA^\top \euPL
\right)
\frac{\partial \euH_{n_o}}{\partial x}(x)
\displaystyle 
\;+\; 
\frac{\partial \euH_{n_o}}{\partial x}(x)^\top \euPL\,  \euB\,   
\frac{\partial L _f^{n_o}h}{\partial x}(x)
+
 \frac{\partial L _f^{n_o}h}{\partial x}(x)^\top \euB ^\top \euPL \frac{\partial \euH_{n_o}}{\partial x}(x)
$\\[0.7em]
from where it follows that
\begin{subequations}\nonumber
\begin{align}
\mathcal{L}_fP(x)
& \leq 
\frac{\partial \euH_{n_o}}{\partial x}(x)^\top\left(
\euPL \euKL \euC + \euC^\top \euKL ^\top\euPL  
- 2q\,  I -\frac{1}{q \nu^2} \euPL \euB \euB ^\top \euPL
\right)\frac{\partial \euH_{n_o}}{\partial x}(x) 
\;+\; 
\frac{\partial \euH_{n_o}}{\partial x}(x)^\top \euPL \euB 
\frac{\partial L _f^{n_o}h}{\partial x}(x) \\
& \hspace{1.2in}
+
 \frac{\partial L _f^{n_o}h}{\partial x}(x)^\top \euB ^\top \euPL \frac{\partial \euH_{n_o}}{\partial x}(x) \\
&\! \leq \!
\frac{\partial \euH_{n_o}}{\partial x}(x)^\top
\euPL \euKL \frac{\partial h}{\partial x}(x)
+ \frac{\partial h}{\partial x}(x)^\top \euKL ^\top\euPL \frac{\partial \euH_{n_o}}{\partial x}(x) 
\\
& \hspace{1in}
 \;-\; q
\left(  2\frac{\partial \euH_{n_o}}{\partial x}(x)^\top \frac{\partial \euH_{n_o}}{\partial x}(x)
-
\nu^2
\frac{\partial L _f^{n_o}h}{\partial x}(x)^\top 
 \frac{\partial L _f^{n_o}h}{\partial x}(x)
\right).
\end{align}
\end{subequations}
Then, using 
\eqref{13},
we get
\begin{eqnarray*}
v^\top \mathcal{L}_fP(x) v\, \leq \,  -q
v^\top\frac{\partial \euH_{n_o}}{\partial x}(x)^\top
\frac{\partial \euH_{n_o}}{\partial x}(x) v 
\leq
-
\frac{q \underline{h}}{\lambda_{\max}(\euPL) \bar{h}}\,
  v^\top P(x) v
\end{eqnarray*}
for all $(x,v)$ such that $\frac{\partial h}{\partial x}(x)v=0$,
which is (\ref{eqn:LfPtangentially}) in the given coordinates).
\end{proof}
\begin{example}
With the above, we see that a Riemannian metric, appropriate for the 
design of
an observer for
the harmonic oscillator with unknown frequency in Example~\ref{ex:RevisitHarmonicOscillator-ComputationofP}, can be 
parameterized on  
$\left(\RR^2\times\RR_{>0}\right)\setminus\left(\{(0,0)\}\times\RR_{>0}\right)$
as
$$
P(x)\!=\! \left(\begin{array}{cccc}
1 & 0 & -x_3 & 0 
\\
0 & 1 & 0 & -x_3
\\
0 & 0 & -x_1& -x_2
\end{array}\right)
\euP
\left(\begin{array}{ccc}
1 & 0 & 0 
\\
0 & 1 & 0 
\\
-x_3 & 0 & -x_1
\\
0 & -x_3 & -x_2
\end{array}\right),
$$
where $\euP $ remains to be designed as a positive 
definite symmetric $4\times 4$ matrix.\hfill $\triangle$
\end{example}

\section{Design of Riemannian Metric $P$ for Lagrangian Systems}
\label{ex:Lagrangian}
In this section, we show
that, besides 
differentially observable systems studied above
Lagrangian systems make another family for which we can easily get an 
expression for a Riemannian metric
that satisfies the differential detectability property introduced in 
Definition~\ref{def:DifferentialDetectability},
at least with symbolic computations and with no need to solve any equation.
To show this,
we follow the ideas in the seminal contribution \cite{Aghannan.Rouchon.03}
and employ the metric used in \cite{Sasaki,Sakai.96}.

Let $\mathcal{Q}$ 
be an $\overline{n}$-dimensional configuration manifold 
equipped with a Riemannian metric $g$. Once we have a chart for 
$\mathcal{Q}$ with coordinates
$q_k$, with $k\in \{1,2,\ldots,\overline{n}\}$, we have also 
coordinates $(q_k,v_l)$ with $(k,l)\in \{1,2,\ldots,\overline{n}\}^2$ 
for its
tangent bundle 
with $q$ being the generalized position and $v$
the generalized velocity. Assume we have a Lagrangian 
$
\mathfrak{L}
:\mathcal{TQ}\to \RR$ of the form
$
\mathfrak{L}(q,v)\;=\; \frac{1}{2}\,   v^\top g(q)\,  v
\;-\;  U(q),$
where the scalar function $U$ is the potential energy. The 
corresponding Euler-Lagrange equations written via any 
chart are
\begin{equation}
\label{LP49}
\dot q _k\;=\; v_k
\quad ,\qquad 
\dot v_l  \;=\; -\mathfrak{C} _{\dummytwo \dummyone }^l  v_\dummytwo 
v_\dummyone \;+\; S_l  (q,t)
\end{equation}
where $k,\,  l  \, ,
\dummytwo\, ,\dummyone \,  
 \in\, \{  1,2,\ldots,\overline{n}\}$;
$S$ is a source term, a known time-varying vector field on
$\reals^{\overline{n}}$;
$\dummytwo $, $\dummyone $ are dummy indices
used for summation in Einstein notation\footnote{
$\sum_m a_m b_{mk}$ is denoted $a_m b_{mk}$ where the fact that the 
index $m$ is used twice means that we should sum in $m$.
}; and
$\mathfrak{C}_{\dummytwo \dummyone }^l$ are the Christoffel symbols associated 
with the metric $g$, namely
$$
\mathfrak{C}_{\dummytwo \dummyone }^l (q)
\;=\; 
\frac{1}{2} 
\left(g(q)^{-1}\right)_{l m} \displaystyle
\left(
\frac{\partial g_{m\dummytwo }}{\partial x_\dummyone }(q) 
+
\frac{\partial g_{m\dummyone }}{\partial x_\dummytwo }(q)
-
\frac{\partial g_{\dummytwo \dummyone }}{\partial x_m}(q)
\right)
\  .
$$
We
consider the measurement $y$ is $q$, namely
$
y \;=\; h(q,v) \;=\;  q.
$

The metric we propose below is for the tangent bundle
$\mathcal{TQ}$.
There are many ways 
of defining
 a Riemannian 
metric for the tangent bundle of a Riemannian manifold \cite{Kowalski.Sekizawa.88}. We follow the same 
route  as the one proposed 
in
 \cite{Aghannan.Rouchon.03} 
to study the local convergence of an observer by considering
the following modification of the Sasaki metric (see \cite[(3.5)]{Sasaki} or \cite[page 55]{Sakai.96}):
$$
P(q,v)\;=\; \left(\begin{array}{cc}
P_{qq}(q,v)& P_{qv }(q,v)
\\
P_{vq}(q,v) & P_{vv}(q,v)
\end{array}\right) 
\  ,
$$
where the  entries of the 
$\overline{n}\times\overline{n}$-dimensional blocks
$P_{qq}$, $P_{qv }$,  $P_{vq }$, and $P_{vv }$ are,
respectively, $P_{ij}$, $P_{i\beta }$, $P_{\alpha j}$, and $P_{\alpha \beta }$,
defined as
\begin{eqnarray*}
P_{ij}(q,v) \!\!\!&=\!\!\!&
ag_{ij}(q)
-c\left(
 g_{i{\dummyone }}(q)\mathfrak{C} _{{\dummytwo }j}^{\dummyone }(q) v_{\dummytwo } 
+
g_{{\dummytwo }j}(q)\mathfrak{C} _{{\dummyone }i}^{\dummytwo }(q) v_{\dummyone } 
\right)\\
& & 
+ bg_{{\dummythree }{\dummyfour }}(q)
\mathfrak{C} _{{\dummytwo }i}^{\dummythree }(q)
\mathfrak{C} _{{\dummyone }j}^{\dummyfour }(q) v_{\dummytwo } v_{\dummyone }
\  ,\\[0.5em]
P_{i\beta }(q,v)&=&
-c g_{i\beta }(q)
+bg_{\beta {\dummyone }}(q)\mathfrak{C} _{{\dummytwo }i}^{\dummyone }(q) v_{\dummytwo }
\  ,
\\[0.5em]
P_{\alpha j}(q,v) & =&
-cg_{\alpha j}(q)
+b g_{\alpha {\dummytwo }}(q)\mathfrak{C} _{{\dummyone }j}^{\dummytwo }(q) v_{\dummyone }
\  ,\\[0.5em]
P_{\alpha \beta }(q,v)&=&
bg_{\alpha \beta }(q)
\  ,
\end{eqnarray*}
where $a$, $b$ and $c$ are strictly positive real numbers satisfying $c^2 < ab$, 
$
g_{\dummytwo \dummyone }
$
are the entries of the metric $g$; and,
here and below,
roman indices $i$, $j$, and $k$ are used to index the components of 
$q$, Greek indices $\alpha$,
$\beta $, and $\gamma $ 
to index the components
of $v$, and $\dummytwo $, $\dummyone $,
$\dummythree $, and $\dummyfour $ are dummy
roman or Greek indices.

We obtain
\\[1em]$
\begin{array}{@{}l@{}}
\left(\begin{array}{cc}
\eta ^\top & \omega ^\top
\end{array}\right) P \left(\begin{array}{c}
\eta  \\ \omega 
\end{array}\right)
=\, \eta _iP_{ij}\eta _j \,+\, \eta _iP_{i\beta  }\omega _\beta  
\,+\, \omega _\alpha P_{\alpha j}\eta _j
\,+\, \omega _\alpha P_{\alpha \beta }\omega _\beta 
\  ,
\\
\hspace{1.4in}=\, a \eta _i g_{ij} \eta _j
\,+\, b
\left(\omega _\alpha +
\mathfrak{C} _{{\dummytwo }i}^{\alpha  }v_\dummytwo \eta _i
\right)g_{\alpha \beta }
\left(\omega _\beta  + \mathfrak{C} _{{\dummyone }j}^{\beta } v_\dummyone \eta _j\right)
 \,-\,
2
 c \eta _i g_{i\beta }\left(\omega _\beta  + \mathfrak{C} 
_{{\dummytwo }j}^{\beta }v_\dummytwo \eta _j\right)
\  .
\end{array}
$\\[1em]
Since $g$ is positive definite and $c^2 < ab$, we see that
$P$ takes positive definite values.

To check that 
we have
the differential detectability property (\ref{eqn:LfPtangentially}),
 we rewrite \eqref{LP49} in the following compact form:
$$
\dot q\;=\; v\quad ,\qquad \dot v\;=\; f_v(q,v,t), \qquad y = h(q,v) = q
\  .
$$
Since we have
$$
{\frac{\partial \h}{\partial (q,v)}(q,v) \!}^{\top}
\frac{\partial \h}{\partial (q,v)}(q,v) 
\;=\; \left(\begin{array}{cc}
I_{\overline{n}} & 0 \\ 0 & 0
\end{array}\right) \in
\reals^{2\overline{n}\times2\overline{n}}
\  ,
$$
 inequality (\ref{eqn:LfPtangentially}) 
is satisfied if we have, for some strictly positive real 
number $\overline{q}$,
\begin{eqnarray*}
\left(\begin{array}{@{}cc@{}}P_{vq} & P_{vv}
\end{array}\right)
\left(\begin{array}{@{}c@{}}
I \\ \displaystyle \frac{\partial  f_v}{\partial v}
\end{array}\right)
+
\left(\begin{array}{@{}cc}
I & \displaystyle \frac{\partial f_v}{\partial v}^\top
\end{array}\right)\left(\begin{array}{@{}c@{}}P_{qv} \\ P_{vv}
\end{array}\right)
+
\frac{\partial P_{vv}}{\partial q } v 
+
\frac{\partial P_{vv}}{\partial v } f_{v }
\; \leq \; -\overline{q}\,  P_{vv}
\  .
\end{eqnarray*}
With the component-wise expression of $f_v$ in (\ref{LP49}),
the symmetry of $g$, and
using Kronecker's delta to denote the identity entries,
the left-hand side above is nothing but 
\\[1em] \vbox{\noindent
$\displaystyle 
\left[
(-cg_{\alpha \dummythree }
+b g_{\alpha {\dummytwo }}\mathfrak{C} _{{\dummyone }\dummythree }^{\dummytwo 
}v_{\dummyone })
\delta _{\dummythree \beta }
-bg_{\alpha \dummytwo }
(\mathfrak{C} _{\dummyone \beta }^\dummytwo + \mathfrak{C} _{\beta \dummyone }^\dummytwo )
 v_\dummyone \right]
$\hfill \null\\[0.3em]\null\hfill$\displaystyle 
\;+\; 
\left[\delta _{\alpha \dummythree }
(
-c g_{\dummythree \beta  }
+bg_{\beta  \dummytwo }\mathfrak{C} _{ \dummyone \dummythree }^\dummytwo v_\dummyone 
)
-(\mathfrak{C} _{\alpha  \dummyone }^\dummytwo + \mathfrak{C} _{\dummyone 
\alpha  }^\dummytwo ) v_\dummyone 
bg_{\dummytwo  \beta }
\right]$\hfill \null\\\null\hfill
$\displaystyle
\;+\;
b\frac{\partial g_{\alpha \beta }}{\partial q_{\dummyone }}v_\dummyone 
$\hfill \null\\\null\hfill
$
\displaystyle = -2cg_{\alpha \beta }
- b\left[
g_{\alpha \dummytwo }\mathfrak{C} _{\beta \dummyone }^\dummytwo 
+
g_{\beta  \dummytwo
\mathfrak{C} _{\alpha  \dummyone }^\dummytwo }
-
\frac{\partial g_{\alpha \beta }}{\partial q_{\dummyone }}
\right]v_\dummyone 
=
-2 c \,  
g_{\alpha \beta }
\  .
$\quad \null  }\\[1em]
Hence,
(\ref{eqn:LfPtangentially}) 
holds since $b$ and $c$ are strictly positive, and the entries of $P_{vv}$ are $b\, g_{\alpha \beta }$.

\begin{example}
Consider a system with
$
\mathfrak{L}(q,v)\;=\; \frac{1}{2}\exp(-2q) v^2$ for all $q, v \in \reals$
as Lagrangian.
The associated metric and its Christoffel symbols are
$
g(q) = \exp(-2q)$, $\mathfrak{C} = -1
$.
Then, the system dynamics are given by
$
\dot q=v$,
$\dot v = v^2$.
Since the (unique) Christoffel symbol is
$\mathfrak{C} = -1$,
we get
$$
P(q,v)\;=\; \exp(-2q)
\left(\begin{array}{cc}
a+2c  v + 
b
v^2
&
-c-bv
\\
-c-bv
&
b
\end{array}\right)
\  .
$$
\hfill $\triangle$
\end{example}

\section{Conclusion}

We have established
 that
strong differential 
detectability
 is already 
sufficient 
for the observer proposed in \cite{Sanfelice.Praly.11}
to guarantee that, at least locally, a
 Riemannian distance between the estimated state and the system state decreases along  solutions.
Moreover in such a case, the existence of a full order
observer implies  the existence of a reduced order one.
This extends the result in \cite[Corollary 3.1]{Besancon.00.SCL} 
established for the 
particular case of an Euclidean metric.

The design of the metric,
exhibiting the strong differential 
detectability
 property
and consequently allowing us to design an observer,
is possible when the system is strongly infinitesimally observable
(i.e., each time-varying linear system
resulting from the linearization along a
solution
 to the system satisfies a uniform observability property).
In such a case, one needs the solution of
an ``algebraic'' (actually a partial differential equation) Riccati 
equation.  This leads to an observer which resembles 
an Extended Kalman Filter.

With the same strong infinitesimal observability property, we can 
also proceed with a linear equation instead of the quadratic Riccati 
equation. In this case the metric we obtain is nothing but an
exponentially weighted observability Grammian. 

The two designs above need the solution of a partial differential 
equation. But thanks to the method of characteristics, it can be 
obtained off-line by  solving ordinary differential equations on a 
sufficiently large time interval and over a grid of initial 
conditions in the system state space. 

A simpler design is possible when the system is
strongly differentially observable (i.e. the mapping state to output 
derivatives is an injective immersion) . Indeed in this case the 
metric can be expressed as a linear combination of functions which 
can be obtained by symbolic computations.
It then remains to choose the linear coefficients.

As already shown in \cite{Aghannan.Rouchon.03}, another case where 
the metric can be obtained via symbolic computations is for Euler-Lagrange 
systems whose Lagrangian is quadratic in the generalized velocities.

Unfortunately, to obtain observers for which convergence holds globally or 
at least regionally and not only locally,
the metric may need to satisfy an extra property. As shown in 
\cite{Sanfelice.Praly.11}, such a property 
can be
a geodesic convexity 
of the level sets of the output function. 
This condition  leads to additional
algebraic equations involving the
Hessian of the output function.

\balance
\bibliographystyle{unsrt} 
\bibliography{long,ObserverGeodesic,RGS}

\appendix

\subsection{Notations and Short glossary of Riemannian geometry}
\label{sec:Glossary}

\begin{enumerate}
\parskip 0pt plus 0pt minus 0pt%
\topsep 0pt plus 0pt minus 0pt
\parsep 0pt plus 0pt minus 0pt%
\partopsep 0pt plus 0pt minus 0pt%
\itemsep 0pt plus 0pt minus 0pt
\settowidth{\labelwidth}{$\bullet$}%
\setlength{\labelsep}{0.5em}%
\setlength{\leftmargin}{\labelwidth}%
\addtolength{\leftmargin}{\labelsep}%
\item
$\SS^n$ denotes the $n$-dimensional unit sphere.
\item
\label{LP65}
Given a function $h:\RR^n\to \RR^p$,
$dh$ denotes its differential form whose expression in coordinates 
$x$ is
$\frac{\partial h_k}{\partial x_j}(x)$
 for each $k$ in $\{1,\ldots,p\}$ and each $j$ in $\{1,\ldots, n\}$.
With $\otimes$, a tensor product, $d\h (x)\otimes d\h(x)$ is a symmetric covariant 
$2$-tensor whose expression in coordinates 
$x$ is $\sum_{k=1}^p \frac{\partial h_k}{\partial x_j}(x) 
\frac{\partial h_k}{\partial x_j}(x)$.
\item\label{item:RiemannianMetric}
A Riemannian metric is a  symmetric covariant $2$-tensor with positive 
definite values.
The associated Christoffel symbols in coordinates $x$ are
$$
\Gamma_{ij}^l\!=\! \frac{1}{2}\sum_{k}
(P^{-1})_{kl}
\left[\frac{\partial P_{ik}}{\partial x_j}
+
\frac{\partial P_{jk}}{\partial x_i}
-
\frac{\partial P_{ij}}{\partial x_k}
\right]
\  .
$$
\item
\label{LP66}
Given a Riemannian metric $P$ and a real valued function $h$, $\gr_P h$ denotes the 
(Riemannian) gradient of $h$.
It is its first covariant derivative. 
Its expression in coordinates $x$ is (see \cite[Sections 1.2 and 
2]{Petersen.06})
$$
\gr_P h(x)\;=\; P(x)^{-1}\frac{\partial h}{\partial x}(x)^\top
\  .
$$

\item
\label{item:CovariantDerivativeOfh}
Given a Riemannian metric $P$ and
a real valued function
 $h$, 
$\Hess _P h$ denotes the (Riemannian) Hessian of $h$. It is its second covariant derivative.
Its expression in coordinates $x$ is
$$
[\Hess _P h (x)]_{ij}=\frac{\partial ^2h}{\partial x_i\partial x_j}(x)
-
\sum_{l}\Gamma _{ij}^l(x)
\frac{\partial h}{\partial x_l}(x)
\  .
$$
It satisfies (see \cite[Sections 1.2 and 2]{Petersen.06})
\begin{equation}
\label{LP96}
\mathcal{L}_{\gr _P h} P(x)\;=\; 2 \,  \Hess _P h (x) 
\  .
\end{equation}
\item
The length of a $C^1$ path
$\gamma $ between points $x_a$ and $x_b$ 
is defined as
$$
\left.\vrule height 1em depth 0.5em width 0pt
L(\gamma )\right|_{s_a}^{s_{b}}\;=\; \int_{s_a}^{s_b}
\sqrt{
\frac{d\gamma }{ds}(s) ^{\top}P(\gamma (s)) \frac{d\gamma }{ds}(s)
} \,  ds,
$$
where
$
\gamma (s_a)\;=\; x_a$ and $\gamma (s_b)\;=\; x_b$.
\item
The Riemannian distance $d(x_a,x_b)$ is 
the minimum of $\left.\vrule height 1em depth 0.5em width 0pt
L(\gamma )\right|_{s_a}^{s_b}$ among all possible piecewise $C^1$ 
paths $\gamma $ between $x_a$ and $x_b$.
A minimizer giving the 
distance is called a
minimizing geodesic and is denoted $\gamma ^*$. 
\item
A topological space equipped with a Riemannian distance is complete
when every geodesic
can be maximally extended to $\RR$.
\item
\label{LP73}
A subset $S$ of $\RR^n$ is said to be weakly geodesically convex
if, for any pair  of points $(x_a,x_b) $ in $S\times S
$, there exists a minimizing geodesic 
$\gamma ^*$
between $x_a=\gamma ^*(s_a)$ and $x_b=\gamma ^*(s_b)$
satisfying
$
\gamma ^*(s)\in S$
for all $s\in [s_a,s_b]$.
A trivial consequence is that any two points in a weakly geodesically convex
can be linked by a minimizing geodesic.
\item
Given a $C^1$ function $h:\RR ^n\mapsto \RR ^p$
and a closed subset $\mathcal{C}$ of $\RR ^n$,
the set
$$
S\;=\; \{x\in \RR ^n:\,  h(x)=0\} \cap \mathcal{C}
$$
is said to be totally geodesic if, for any pair $(x,v)$ in 
$S\times\RR  ^n$ such that
$
\frac{\partial h}{\partial x}(x)\,  v\;=\; 0$ and
$v^\top P(x)\,  v\;=\; 1$,
any geodesic $\gamma $ with
$
\gamma (0)\;=\; x$,
$\frac{d\gamma }{ds}(0)\;=\; v
$
satisfies
$
h(\gamma (s)) \;=\; 0$
for all $s \in J_\gamma$,
where $J_\gamma $ is the maximal interval containing $0$ so that $\gamma (J_\gamma )$ 
is contained in $\mathcal{C}$.
\item
\label{LP63}
Given a set of coordinates for $x$, the Lie derivative $\mathcal{L}_fP$
of a symmetric covariant $2$-tensor 
$P$ is, for all $v$ in $\RR ^n$,
\begin{eqnarray*}
\displaystyle v^\top \mathcal{\mathcal{L}}_f P(x)\,  v &  = &
\lim_{t\to 0} 
\left[\frac{[(I+t\frac{\partial f}{\partial x}(x))v]^\top
P(X(x,t)) 
[(I+t\frac{\partial f}{\partial x}(x))v]}{t}- \frac{v^\top P(x) v}{t}\right]\\
&=&
\frac{\partial}{ \partial \x}\left(
\vrule height 1em depth 0.5em width 0pt
v^{\top} P(\x) \,  v\right)\,  \f(x) 
\;+\; 2\,  
v^{\top} P(\x)\left(\frac{\partial \f}{\partial \x}(\x) \,  v\right)
\end{eqnarray*}
where $t\mapsto X(x,t)$ is the solution to \eqref{eqn:Plant1}.
If there exist coordinates in $\RR^n$ denoted $x$ and a function $\varphi :\RR^n\to 
\RR^p$ such that the expression of $P$ is
$$
P(x)\;=\; \frac{\partial \varphi }{\partial x}(x)^\top \euP \frac{\partial \varphi }{\partial x}(x)
$$
where $\euP$ is a symmetric matrix, then we have
\begin{equation}
\label{LP64}
\mathcal{L}_fP(x)\;=\; 
\frac{\partial L_f\varphi }{\partial x}(x)^\top \euP \,  \frac{\partial \varphi }{\partial x}(x)
\;+\; 
\frac{\partial \varphi }{\partial x}(x)^\top \euP \,  \frac{\partial L_f\varphi }{\partial x}(x)\ ,
\end{equation}
where $L_f\varphi $ is the image by $\varphi $ of the vector field $f$ 
(in $\RR^n$). Indeed, we have
\begin{eqnarray*}
v^\top \mathcal{L}_fP(x)\,  v&
=&\displaystyle 
2\,  v^\top\frac{\partial \varphi }{\partial x}(x)^\top \euP\,  
\frac{\partial L_f\varphi }{\partial x}(x)v\ .
\end{eqnarray*}
We would like the reader to distinguish the notation
$\mathcal{L}_fP$ for the Lie derivative of a symmetric covariant $2$-tensor from
$L_f\varphi $, which is used for the more usual Lie derivative of a function $\varphi $, or 
equivalently, the vector field induced by a function.
\end{enumerate}

\end{document}